\documentclass[11pt,a4paper]{article}
\usepackage[colorlinks=true,linkcolor=blue]{hyperref}
\usepackage{cmap}
\usepackage[utf8x]{inputenc}
\usepackage[T1]{fontenc}
\usepackage[english]{babel}
\usepackage{amssymb,amsmath}
\usepackage{amsthm}
\usepackage{bbm}
\usepackage{microtype}
\usepackage{lmodern}

{
\rmfamily
\DeclareFontShape{T1}{lmr}{m}{sc}{<->ssub*cmr/m/sc}{}
\DeclareFontShape{T1}{lmr}{b}{sc}{<->ssub*cmr/b/sc}{}
\DeclareFontShape{T1}{lmr}{bx}{sc}{<->ssub*cmr/bx/sc}{}
}

\newcommand{\thmheadercommand}[1]{\textbf{\scshape{}#1.\\*}}

\newtheoremstyle{yannthm}{\topsep}{\topsep}{\slshape}{}{\scshape\bfseries}{.}{.5em}{%
\thmname{#1}\thmnumber{ #2}\thmnote{#3}%
}
\newtheoremstyle{yannthm2}{\topsep}{\topsep}{}{}{\scshape\bfseries}{.}{.5em}{%
\thmname{#1}\thmnumber{ #2}\thmnote{#3}%
}

\def\d{\operatorname{d}\!{}}

\def\R{{\mathbb{R}}}

\renewcommand{\geq}{\geqslant}

\newcommand{\deq}{\mathrel{\mathop{:}}=}

\newcommand{\from}{\colon} % correct ':' in f\from X \to Y
\renewcommand{\epsilon}{\varepsilon}
\renewcommand{\phi}{\varphi}

\DeclareMathOperator{\Cov}{Cov}
\let\oldPr\Pr
\renewcommand{\Pr}{\oldPr\nolimits}
\newcommand{\E}{\mathbb{E}}

\DeclareMathOperator{\Id}{Id}

\newcommand{\norm}[1]{\left\lVert#1\right\rVert}

\newenvironment{dem}[1][]{\begin{proof}[\thmheadercommand{Proof#1}]~\newline\ignorespaces}{\end{proof}}

{
\theoremstyle{yannthm}
\newtheorem{defi}{Definition}
\newtheorem*{defi*}{Definition}
\newtheorem{prop}[defi]{Proposition}
\newtheorem*{prop*}{Proposition}
\newtheorem{thm}[defi]{Theorem}
\newtheorem*{thm*}{Theorem}
\newtheorem{lem}[defi]{Lemma}
\newtheorem*{lem*}{Lemma}
\newtheorem{cor}[defi]{Corollary}
\newtheorem*{cor*}{Corollary}

\newtheorem*{ex*}{Example}

\newtheorem*{subenonce}{}

\theoremstyle{yannthm2}

\newtheorem*{exo*}{Exercise}

\newtheorem*{rem*}{Remark}

\newtheorem*{subenonce2}{}

}

\newenvironment{enonce}[1]{\begin{subenonce}[#1]}{\end{subenonce}}

\newcommand{\transp}[1]{#1^{\!\top}\!}

\title{The Extended Kalman Filter is a Natural Gradient Descent in
Trajectory Space}

\author{Yann Ollivier}
\date{}

\newcommand{\gaussian}{\mathcal{N}}

\newcommand{\deltatheta}{\delta\hspace{-.05em}\theta}

\newcommand{\traj}{\mathbf{s}}
\newcommand{\Traj}{\mathcal{S}}
\newcommand{\obs}{\mathbf{y}}
\newcommand{\pobs}{p_{\mathrm{obs}}}

\newcommand{\Jt}[2]{J_{#1\downarrow #2}}

\begin{document}

\maketitle

\begin{abstract}
The extended Kalman filter is perhaps the most standard tool to estimate
in real time the state of a dynamical system from noisy measurements of
some function of the system, with extensive practical applications (such as position tracking via GPS).
While the plain Kalman filter for linear systems is well-understood, the
extended Kalman filter relies on linearizations which have been debated.

We recover the exact extended Kalman filter equations from first
principles in statistical learning: the extended Kalman filter is equal
to Amari's \emph{online natural gradient}, applied in the space of
trajectories of the system.  Namely, each possible trajectory of the
dynamical system defines a probability law over possible observations. In
principle this makes it possible to treat the underlying trajectory as
the parameter of a statistical model of the observations. Then the
parameter can be learned by gradient ascent on the
log-likelihood of observations, as they become available.  Using Amari's
\emph{natural} gradient from information geometry (a gradient
descent preconditioned with the Fisher matrix, which provides
parameterization-invariance) exactly recovers the extended Kalman filter.

This applies only to a particular choice of process noise in the Kalman
filter, namely, taking noise proportional to the posterior covariance—a
canonical
choice in the absence of specific model information.
\end{abstract}

\paragraph{Overview.}
State estimation consists in estimating the current state of a dynamical
system given noisy observations of a function of this system. Namely,
consider a dynamical system
with state $s_{t}$, inputs $u_{t}$ and dynamics $f$, namely,
\begin{equation}
\label{eq:syst}
s_{t}=f(s_{t-1},u_{t})
\end{equation}
and assume we have access to noisy observations $y_t$ of some function
$h$ of the system,
\begin{equation}
\label{eq:obs}
y_t=h(s_{t},u_{t})+\gaussian(0,R)
\end{equation}
with covariance matrix $R$. One of the main problems of filtering theory
is to estimate the current state $s_t$ given the observations $y_t$
(assuming that $f$, $h$, $R$, and the inputs or control variables $u_t$ are
known).

We prove the exact equivalence of two methods to tackle this problem:
\begin{itemize}
\item The \emph{extended Kalman filter}, the most standard tool designed
to deal with this problem: it is built in a Bayesian setting as a
real-time
approximation of the posterior mean and covariance of the state given the
observations. (We use a particular variant of the filter
where the process noise on $s_t$ is modeled as proportional to the
posterior covariance, $Q_t\propto P_{t|t-1}$ in Def.~\ref{def:kalman}
\cite[\S3.2.2]{Nelson_thesis} \cite[\S5.2.2]{Haykin_book},
a canonical choice in the absence of further information. This choice
introduces ``fading memory'' which robustifies the filter
\cite[\S5.5]{simon2006kalmanbook}.)
\item The \emph{online natural gradient}, a classical tool from statistical
learning to estimate the parameters of a probabilistic model. Here, the
hidden parameter to be estimated is the whole \emph{trajectory}
$\traj=(s_t)_{t\geq 0}$. Letting $\Traj$ be the set of trajectories of
\eqref{eq:syst}, each possible trajectory $\traj=(s_t)\in \Traj$ defines a
probability distribution $p(\obs|\traj)$ on observation sequences
$\obs=(y_t)_{t\geq 1}$ via the observation model \eqref{eq:obs}. So $\traj$ can be seen as the
parameter of a probabilistic model on $\obs$. Then, in principle, $\traj$
can be
learned by online gradient descent $\frac{\ln p(y_t|\traj)}{\partial \traj}$ in
the space of trajectories:
each time a new observation $y_t$ becomes available, one can re-estimate
$\traj$ using a gradient step on the log-likelihood of $y_t$
knowing $\traj$.

The \emph{natural} gradient descent \cite{Amari1998} preconditions the
gradient steps
by the inverse Fisher matrix of the model
$p(\obs|\traj)$ with respect to $\traj$. This is motivated by invariance
to changes of variables over which the model is expressed, and by
theorems of asymptotic optimality \cite{Amari1998}.
\end{itemize}

We claim that these two methods yield the same estimate of $s_t$ at time
$t$ (Thm.~\ref{thm:kalnat}). The same holds in continuous time for the extended Kalman--Bucy
filter (Thm.~\ref{thm:natkalbucy}).

This largely extends a previous result by the author, which dealt
with the case $f=\Id$: namely, it was shown in \cite{natkal} that the natural
gradient descent to estimate the parameter $\theta$ of a probabilistic
model from observations $y_t$, is equivalent to applying a Kalman filter
to the hidden state $s_t=\theta$ for all $t$. Thus the previous result
viewed the natural gradient as a particular case of an extended Kalman
filter with ``static'' dynamics; here we view the extended Kalman filter
as a natural gradient descent in the space of trajectories, and recover
the previous result when $f=\Id$.

\bigskip

This result may contribute to the understanding of the extended Kalman
filter. The use of Kalman-like filters in navigation systems (GPS,
vehicle control, spacecraft...), time series analysis, econometrics,
etc.\ \cite{sarkka_book}, is extensive to the point it has been described
as one of the greater discoveries of mathematical engineering
\cite{kalmanmatlab}.  But while the plain Kalman filter (which deals with
linear $f$) is exactly optimal, the extended Kalman filter relies on
linear expansions. Variants of the extended filter have been proposed,
for instance using higher-order expansions for certain terms, though
with more limited use \cite{rg2011secondorderkalman}.
On the other hand,
the natural gradient can be constructed from first principles.
Remarkably, the quite complicated formulas defining the extended Kalman
filter can be derived exactly from its natural gradient interpretation.

\bigskip

However, two technical points make the precise statement of the
correspondence (Theorem~\ref{thm:kalnat}) more subtle.

First, an important choice when applying the extended Kalman filter is
the choice of system noise $\gaussian(0,Q)$ that is added to the dynamical
system \eqref{eq:syst}.  (One may think of the process noise $Q$ in the
Kalman filter either as actual noise in a stochastic system, or as a
modeling tool to apply the Kalman filter when knowledge of the
deterministic system $f$ is imperfect; the results below hold regardless of 
interpretation.) Often, $Q$ is adjusted by trial and error. A canonical
choice is to take $Q$ proportional to
the posterior covariance on $s$ (Def.~\ref{def:Q}, \cite[\S3.2.2]{Nelson_thesis} \cite[\S5.2.2]{Haykin_book}); this is equivalent to introducing
\emph{fading memory} into the filter \cite[\S7.4]{simon2006kalmanbook}.
% Fagin 1964 "recursive
% linear regression theory", Morrison 1969 "introduction to sequential
% smoothing and prediction"

Our result applies only in
the latter case; this is certainly a restriction. Fundamentally, choices
such as $Q=\Id$
define a preferred basis in state space, while the extended Kalman
filter with $Q$ proportional to the posterior variance can be expressed
in an abstract, basis-free vector space. Since the natural gradient is
basis-invariant, it can only be equivalent to another basis-invariant
algorithm.\footnote{Other choices of $Q$, such as  $Q=\Id$, do have an
interpretation as gradient descents in
trajectory space, but using quite artificial preconditioning matrices instead
of the Fisher matrix; we do not develop this point.}

Our results relate the extended Kalman filter with nonzero $Q$ to the
natural gradient over trajectories of the \emph{noiseless} system
\eqref{eq:syst}. The choice of noise $Q$ for applying the Kalman filter
corresponds to different natural gradient learning rates: the particular
choice $Q=0$ corresponds to a learning rate $1/t$ in the natural
gradient, while positive $Q$ correspond to larger learning rates. 

Second, the natural gradient uses quantities expressed in an abstract
Riemannian manifold of trajectories $\traj$; still, to perform an actual
update of $\traj$, a numerical representation of $\traj$ has to be used.
(The direction of the natural gradient is parameterization-invariant, but
the actual step requires an explicit parameterization, whose influence
vanishes only in the limit of small learning rates.) The space of
trajectories $\traj$ could be parameterized, for instance, by the initial
state $s_0$, or the state $s_t$ at any time $t$ provided $f$ is
invertible. The correspondence turns out to be exact if, when the
observation $y_t$ becomes available at time $t$, the natural gradient
update uses the current state $s_t$ to parameterize of the trajectory
$\traj$. One one hand this seems quite natural, and computationally
convenient at time $t$; on the other hand, it means we are performing a natural
gradient descent in a coordinate system that shifts in time.

% This effect is specific to discrete time and disappears in continuous
% time (Section~\ref{sec:cont}).  Indeed, in continuous time the updates to
% the trajectory are infinitesimal at each instant, and the natural
% gradient is exactly parameterization-invariant. In that case, the online
% natural gradient is directly defined over the manifold of trajectories,
% and provides a derivation of the Kalman--Bucy filter (with process noise $Q_t$
% proportional to $P_t$) from first principles.

\paragraph{Example: recovering the natural gradient from the extended
Kalman filter for statistical learning problems.} The correspondence
works both ways: in particular, it can be used to view the online natural
gradient on a parameter $\theta$ of a statistical model, as a particular
instance of extended Kalman filtering. This important example corresponds
to $f=\Id$ above, and is the case treated in \cite{natkal}. We summarize
it again for convenience.

Let $p(y_t|u_t,\theta)$ be a statistical model to predict a quantity
$y_t$ from an input $u_t$ given a parameter $\theta$. We assume that the
model can be written as $y_t\sim \pobs(y_t|h(\theta,u_t))$ where
$h(\theta,u_t)$ is a function that encodes the prediction on $y_t$, and
the noise model $\pobs$ is an exponential family with mean parameter
$h(\theta,u_t)$, such as $y_t=h(\theta,u_t)+\gaussian(0,R)$. The function
$h$ may be anything, such as $h(\theta,u_t)=\transp{\theta} u_t$ for a
linear model, or a feedforward neural network with input $u_t$ and
parameters $\theta$.

A standard approach for this problem would be stochastic gradient
descent: updating the parameter $\theta$ via gradient descent of $\ln
p(y_t|u_t\theta)$ for each new observation pair $(u_t,y_t)$.
But the extended Kalman filter can also be applied to this
problem by viewing
$\theta$ as the hidden state of a static system, namely, $s_t=\theta$ and
$f=\Id$, and treating the $y_t$ as observations of $\theta$ knowing
$u_t$. See eg \cite{singhalwu1988} for an early example with neural networks.
Following \cite{natkal}, we extend the extended Kalman
filter in Def.~\ref{def:kalman} to cover any exponential family as the model for $y_t$ given
$h(s_t,u_t)$: this allows the Kalman filter to deal with
discrete/categorical data $y_t$, for instance, by letting $h(\theta,u_t)$
be the list of probabilities of all classes.

The main result from \cite{natkal} states that the extended Kalman filter
for this problem, is exactly equivalent to the online natural gradient on
$\theta$. This is a corollary of the present work by taking $f=\Id$ and
$s_t=\theta$: indeed, with $f=\Id$ we can identify the set of
trajectories $\traj\in \Traj$ with their value at any time, and the gradient
descent on $\traj$ becomes a gradient on $\theta$.

So the online natural gradient for a statistical problem with parameter
$\theta$ appears as a particular instance of the extended Kalman filter
on a static system $f=\Id$, while the extended Kalman filter for general
$f$ appears as a particular case of the online natural gradient in the
more abstract space of trajectories.

\paragraph{Does this provide a convergence proof for the extended Kalman
filter, via the theory of stochastic gradient descent?} Not really, as consecutive observations in a dynamical system
are not independent and identically distibuted. The online natural
gradient on a dynamical system is not quite an instance of stochastic
gradient descent.

\paragraph{Related work.} The role of the information matrix in Kalman
filtering was recognized early \cite[\S7.5]{jazwinski_book}, and led to
the formulation of the Kalman filter using the inverse covariance matrix
known as the ``information filter'' \cite[\S6.2]{simon2006kalmanbook}.
However, except in the static case treated in \cite{natkal}, this does
not immediately translate into an equivalence between extended Kalman
filtering and natural gradient descent, as is clear from the amount of
work needed to prove our results.

Several recent works make a link between Kalman filtering and preconditioned
gradient descent in some particular cases.
% \cite{singhalwu1988} is an early example of the use of Kalman filtering
% for training feedforward neural networks in statistical learning, but
% does not mention the natural gradient.
\cite{ruck1992comparative} argue
that for neural networks, backpropagation, i.e., ordinary gradient
descent, ``is a degenerate form of the extended Kalman filter''.
\cite{bertsekas96} identifies the extended Kalman filter with a
Gauss--Newton gradient descent for the specific case of nonlinear
regression.  \cite{freitas2000hierarchical} interprets process noise in
the static Kalman filter as an adaptive, per-parameter learning rate,
thus akin to a preconditioning matrix.  \cite{simandl2001CramerRao} uses
the Fisher information matrix to study the variance of parameter
estimation in Kalman-like filters, without using a natural gradient;
\cite{bottoulecun2003} comment on the similarity between Kalman filtering
and a version of Amari's natural gradient for the specific case of least
squares regression; \cite{martensnatgrad} and \cite{gradnn} mention the
relationship between natural gradient and the Gauss--Newton Hessian
approximation; \cite{patel_kalmansgd} exploits the relationship between
second-order gradient descent and Kalman filtering in specific cases
including linear regression; \cite{li2017information} use a natural
gradient descent over Gaussian distributions for an auxiliary problem
arising in Kalman-like Bayesian filtering, a problem independent from the
one treated here.

\cite{freshkalman} interpret the Kalman filter as an online Newton method
over a variable representing the trajectory. Namely, defining the ``past
trajectory'' of the system as $z_t\deq (s_1,\ldots,s_t)$, and denoting
the log-likelihood function by
$J_t(z_t)\deq \frac12 \sum_{s=1}^t
\norm{s_t-f(s_{t-1},u_t)}^2_{Q^{-1}}+\frac12 \sum_{s=1}^t
\norm{y_t-h(s_t,u_t)}^2_{R^{-1}}$, they prove that the Kalman filter can
be seen, at each time step, as one step of the Newton method on $z_t$ to find the minimum of $J_t$.
This is somewhat reminiscent of the approach taken here.  However, the
derivation for the nonlinear case is incomplete
% \footnote{Indeed, it is
% assumed in \cite{freshkalman} that the trajectory estimate $\hat z_{t-1}$ at time $t-1$ is the
% minimizer of $J_{t-1}$, from which it is proved that the Kalman filter at
% step $t$ is equivalent to a Newton step for $\hat z_t$; but this does not
% imply that $\hat z_t$ is a minimizer of $J_t$ unless $J_t$ is quadratic,
% so the induction is incomplete except in the linear case.}
(otherwise,
this
would prove optimality of the extended Kalman filter even in the
nonlinear case). Their result states that assuming $\hat z_{t-1}$
minimizes $J_{t-1}$, then the extended Kalman filter tries to find the state $s_t$
that minimizes $J_t$, via
one Newton step. In the linear case, $J_t$ is quadratic, so this Newton step successfully
finds the minimum of $J_t$, so $\hat z_t$ minimizes $J_t$ and the idea can be iterated. Therefore the plain
(non-extended) Kalman filter can be seen as an online Newton method on
$z_t=(s_1,\ldots,s_t)$. However, if $\hat z_t$ does not exactly minimize $J_t$
then the extended Kalman filter at time $t+1$ does not coincide with a
Newton step anymore. Using the Fisher information matrix instead of the
Hessian, namely, a natural gradient instead of a Newton method, helps with this issue.

\paragraph*{Notation conventions.} In
statistical learning, the external inputs or regressor variables are often
denoted $x$%, and we will apply a Kalman filter to a system whose state is
%$\theta$, the parameter to be learned
. In Kalman filtering, $x$ often denotes the state of
the system, while the external inputs are often $u$. Thus %, to avoid conflict of notation, 
we will avoid $x$ altogether and denote by $u$ the inputs and
by $s$ the state of
the system.

The variable to be predicted at time $t$ will be $y_t$, and $\hat y_t$ is
the corresponding prediction. In general $\hat y_t$ and $y_t$ may be
different objects in that $\hat y_t$ encodes a full probabilistic
prediction for $y_t$. For Gaussians with known variance, $\hat y_t$ is
just the predicted mean of $y_t$, so in this case $y_t$ and $\hat y_t$
are the same type of object.  For Gaussians with unknown variance, $\hat
y$ encodes both the mean and second moment of $y$. For discrete
categorical data, $\hat y$ encodes the probability of each possible
outcome $y$.

% Thus, the formal setting for this text is as follows: we are given a
% sequence of finite-dimensional observations $(y_t)$ with each $y_t\in
% \R^{\dim(y)}$, a sequence of inputs $(u_t)$ with each $u_t\in
% \R^{\dim(u)}$, a parametric model $\hat y=h(\theta,u_t)$ with parameter
% $\theta\in \R^{\dim(\theta)}$ and $h$ some fixed smooth function from
% $\R^{\dim(\theta)}\times \R^{\dim(u)}$ to $\R^{\dim(\hat y)}$. We are
% given an exponential family (output noise model) $p(y|\hat y)$ on $y$
% with mean parameter $\hat y$ and sufficient statistics $T(y)$ (see the
% Appendix), and we define the loss function $\ell_t\deq -\ln p(y_t|\hat y_t)$.

The natural gradient descent on parameter $\theta_t$ will use the Fisher
matrix $J_t$. The Kalman filter will have posterior covariance matrix
$P_t$.

For multidimensional quantities $x$ and $y=f(x)$, we denote by $\frac{\partial
y}{\partial x}$ the Jacobian matrix of $y$ w.r.t.\ $x$, whose $(i,j)$
entry is $\frac{\partial f_i(x)}{\partial x_j}$. This satisfies the chain
rule $\frac{\partial z}{\partial y}\frac{\partial y}{\partial
x}=\frac{\partial
z}{\partial x}$. With this convention, 
gradients of real-valued functions are \emph{row} vectors, so that a
gradient descent takes the form $x\gets x - \eta\, \transp{(\partial
f/\partial x)}$.

For a column vector $u$, $u^{\otimes 2}$ is synonymous with
$u\transp{u}$, and with $\transp{u}u$ for a row vector.

%We always implicitly refer to the ``extended'' (nonlinear)
%version of the Kalman filter.

\section{Natural Gradient Descent}
\label{sec:natgrad}

A standard approach to optimize the
parameter $\theta$ of a probabilistic model, given a sequence of
observations $(y_t)$, is an online gradient descent
\begin{equation}
\theta_t\gets \theta_{t-1} + \eta_t \transp{\frac{\partial \ln p(y_t|\theta)}{\partial
\theta}}
\end{equation}
with learning rate $\eta_t$.  This simple gradient descent is
particularly suitable for large datasets and large-dimensional models
\cite{bottoulecun2003}, and has become a staple of current statistical
learning, but has several
practical and theoretical shortcomings. For instance, it uses the same
non-adaptive learning rate for all parameter components. Moreover, simple
changes in parameter encoding or in data
presentation (e.g., encoding black and white in images by 0/1 or 1/0) can
result in different learning performance.

This motivated the introduction of the
\emph{natural gradient} \cite{Amari1998}. 
 It is built to achieve
invariance with respect to parameter re-encoding; in particular,
learning become insensitive to the characteristic scale of each parameter
direction, so that different directions naturally get suitable learning
rates. The natural gradient is the only general way to achieve
such invariance \cite[\S2.4]{Amari2000book}.

The natural gradient preconditions
the gradient descent with
$J(\theta)^{-1}$ where $J$ is the \emph{Fisher information matrix} 
\cite{Kullback} with
respect to the parameter $\theta$. For a smooth probabilistic model
$p(y|\theta)$ over a random variable $y$ with parameter $\theta$, the
latter
is defined as
\begin{equation}
J(\theta)\deq \E_{y\sim p(y|\theta)} \left[{\frac{\partial \ln
p(y|\theta)}{\partial \theta}}^{\otimes 2}\right]=-\E_{y\sim
p(y|\theta)}\left[
\frac{\partial^2\ln p(y|\theta)}{\partial \theta^2}\right]
\end{equation}
If the model for $y$ involves an input $u$, then an additional expectation or
empirical average over the input is introduced in the definition of $J$
\cite[\S8.2]{Amari2000book} \cite[\S5]{martensnatgrad}.

Intuitively, $J$ captures the change in the distribution $p_\theta$ when
$\theta$ changes infinitesimally, measured by the relative entropy
(Kullback--Leibler divergence), namely
\begin{equation}
\mathrm{KL}(p_{\theta+\deltatheta}|p_\theta)=\frac12
\transp{\deltatheta}J(\theta)\deltatheta+O(\deltatheta^3)
\end{equation}
In particular, this only depends on $\theta$ via $p_\theta$. Namely, making a
change of variables in the parameters $\theta$ of a probabilistic model
$p_\theta$ will not change
$\mathrm{KL}(p_{\theta+\deltatheta}|p_\theta)$; the norm 
$\transp{\deltatheta}J(\theta)\deltatheta$ is a parameter-invariant way
to measure the change of $p_\theta$ induced by $\theta$, and this turns
$\Theta$ into a Riemannian manifold \cite{Amari2000book}.

Intuitively the natural gradient is thus the steepest gradient direction
in Kullback--Leibler distance: the natural gradient direction
$J(\theta)^{-1}\partial \ell(\theta)/\partial \theta$ of a function
$\ell(\theta)$ gives, at first order, the direction $\deltatheta$ with
steepest increase of $f$ for the minimum change
$\mathrm{KL}(p_{\theta+\deltatheta}|p_\theta)$ of $p_\theta$ \cite{igo}.

However, this comes at a large computational cost for large-dimensional
models: just storing the Fisher matrix already costs $O((\dim
\theta)^2)$. Various strategies are available to approximate the natural
gradient for complex models such as neural networks, using diagonal or
block-diagonal approximation schemes for the Fisher matrix, e.g.,
\cite{TONGA,gradnn,riemaNN,grosse2015scaling,martens2015optimizing}.

Definition~\ref{def:natgrad} below formally
introduces the \emph{online} natural gradient.

\begin{defi}[ (Online natural gradient)]
\label{def:natgrad}
Consider a statistical model with parameter $\theta$ that predicts an
output $y$ given an input $u$, via a model $y\sim p(y|u,\theta)$.
Given observation pairs $(u_t,y_t)$, the goal is to minimize, online, the
log-likelihood loss function
\begin{equation}
-\sum_t \ln p(y_t|u_t,\theta)
\end{equation}
as a function of $\theta$.

The \emph{online natural gradient} maintains a current estimate
$\theta_t$ of the parameter $\theta$, and a current approximation $J_t$ of the
Fisher matrix. The parameter is estimated by a gradient descent with
preconditioning matrix $J_t^{-1}$, namely
\begin{align}
\label{eq:natgradJ}
J_t &\gets (1-\gamma_t) J_{t-1} + \gamma_t \,\E_{y\sim 
p(y|u_t,\theta)} \left[\frac{\partial \ln p(y|u_t,\theta)}{\partial
\theta}^{\otimes 2}\right]
\\
\label{eq:natgradtheta}
\theta_t &\gets 
% \theta_{t-1} -\eta_t \, J_t^{-1}\transp{\left(
% \frac{\partial \ell_t(\theta)}{\partial \theta}
% \right)}
% =
\theta_{t-1} +\eta_t \, J_t^{-1}\transp{\left(
\frac{\partial \ln p(y_t|u_t,\theta)}{\partial \theta}
\right)}
\end{align}
with learning rate $\eta_t$ and Fisher matrix decay rate $\gamma_t$.
\end{defi}

In the Fisher matrix update, the expectation over all possible values
$y\sim p(y|\hat y)$ can often be computed algebraically (for a given
input $u_t$), but this is
sometimes computationally bothersome (for instance, in neural networks,
it requires $\dim(\hat y_t)$ distinct backpropagation steps \cite{gradnn}). A
common solution \cite{APF00,TONGA,gradnn,PB13} is to just use the value $y=y_t$ (\emph{outer
product} approximation) instead of the expectation over $y$. Another is to use a
Monte Carlo approximation with a single sample of $y\sim p(y|\hat y_t)$
\cite{gradnn,riemaNN}, namely, using the gradient of a synthetic sample
instead of the actual observation $y_t$ in the Fisher matrix. These
latter two solutions are often confused; only the latter provides an
unbiased estimate, see discussion in \cite{gradnn,PB13}.

The online ``smoothed'' update of the Fisher matrix in
\eqref{eq:natgradJ} reuses Fisher matrix values computed at previous
values of $\theta_t$ and $u_t$, instead of using the exact Fisher matrix
at $\theta_t$ in \eqref{eq:natgradtheta}.
Such or similar
updates are used in \cite{TONGA,riemaNN}. The reason is at least
twofold. First, the exact Fisher matrix involves an
expectation over the inputs $u_t$ \cite[\S8.2]{Amari2000book}, so using
it would
mean recomputing the value of the Fisher
matrix on all previous observations each time $\theta_t$ is updated.
Instead, to keep
the algorithm online, \eqref{eq:natgradJ} reuses values computed on
previous observations $u_t$, even though they were computed using an
out-of-date parameter $\theta$. The decay rate $\gamma_t$ controls this
moving average over observations
(e.g., $\gamma_t=1/t$ realizes an
equal-weight average over all inputs seen so far).
%But $\gamma_t=1$ works if model does not depends on any inputs u_t
%(i.e., for unsupervised learning of $y_t$).
Second, the expectation over $y\sim p(y|u_t,\theta)$ in \eqref{eq:natgradJ}
is often replaced with a Monte Carlo estimation with only one value of
$y$, and averaging over time compensates for this Monte Carlo sampling.

As a consequence, since $\theta_t$ changes over time, this means that the
estimate $J_t$ mixes values obtained at different values of $\theta$, and
converges to the Fisher matrix only if $\theta_t$ changes slowly, i.e.,
if $\eta_t\to 0$.  The correspondence below with Kalman filtering
suggests using $\gamma_t=\eta_t$.

\paragraph{Natural gradient descent in different charts of a manifold.}
One motivation for natural gradient is its invariance to a change of
parameterization of the model REF. However, this holds only in the limit
of small learning rates $\eta_t\to 0$, or in continuous time; otherwise
this is true only up to
$O(\eta_t^2)$. Indeed, if $\theta$ belongs to a
manifold, the object $J_t^{-1} \transp{\left(
\frac{\partial \ell_t}{\partial \theta}\right)}$ is a well-defined
tangent vector at $\theta$, but the additive update $\theta\gets
\theta-\eta_t \,J_t^{-1} \transp{\left(
\frac{\partial \ell_t}{\partial \theta}\right)}$ is still performed on an
explicit parameterization (chart).\footnote{A possible solution is to use
the geodesics of the Riemannian manifold defined by the Fisher metric,
but this is rarely convenient except in particular situations where these
geodesics are known explicitly (e.g.,
\cite{bensadon2015gigo,bonnabel2013stochastic}).} Thus, each time an explicit update is
performed, a coordinate system must be chosen.

\newcommand{\chart}{\Phi} 
\newcommand{\Tang}{\mathbf{T}}
\newcommand{\abstheta}{\vartheta}
\newcommand{\numtheta}{\theta}
\newcommand{\absJ}{\mathcal{J}}
\newcommand{\numJ}{J}
\newcommand{\absloss}{\mathcal{L}}
\newcommand{\numloss}{\ell}
\newcommand{\absp}{\mathfrak{p}}
\newcommand{\nump}{p}

Thus, from now on we will explicitly separate the abstract points
$\abstheta\in \Theta$ in the abstract parameter manifold $\Theta$, and their
expression $\numtheta\in \R^{\dim(\Theta)}$ in a coordinate system.
Likewise, we will denote $\numJ$ the Fisher matrix in a coordinate
system, and $\absJ$ the corresponding abstract Fisher metric,
a $(0,2)$-tensor on $\Theta$. We will denote $\absp_t(y_t|\abstheta)$ the
probability distribution on observations at time $t$ knowing the
parameter $\abstheta\in \Theta$, and $\nump_t(y_t|\numtheta)$ the same
model in a coordinate system. The loss function to be minimized is
$-\sum_t \ln
\absp_t(y_t|\abstheta)$.

In particular, when observing $y_t$ at time $t$, the natural gradient
direction with Fisher metric tensor $\absJ$ is
\begin{equation}
\absJ^{-1}\,\frac{\partial \ln \absp_t(y_t|\abstheta)}{\partial\abstheta}
\end{equation}
where $\frac{\partial \ln \absp_t(y_t|\abstheta)}{\partial\abstheta}$ is
a cotangent vector, which becomes a tangent vector after applying
$\absJ^{-1}$.
Then we would like to consider an update of the type
\begin{equation}
\abstheta\gets \abstheta + \eta_t \,\absJ^{-1}\,\frac{\partial \ln
\absp_t(y_t|\abstheta)}{\partial\abstheta}
\end{equation}
with learning rate $\eta_t$. However, this $+$ sign does not make
sense in a manifold, so we have to apply this in an explicit
parameterization (chart), then jump back to the manifold.\footnote{This
only matters at
second order in the learning rate: two different parameterizations will provide
updates differing by $O(\eta_t^2)$. In particular, in continuous time
these considerations disappear, and the continuous-time trajectory
\begin{equation}
\frac{\d \abstheta}{\d t}=\absJ^{-1}\,\frac{\partial \ln
\absp_t(y_t|\abstheta)}{\partial\abstheta}
\end{equation}
is parameterization-independent.
}

This is the object of the next definition: at each step, we first jump to a
chart $\chart$, apply the natural gradient update in that chart, and jump
back to the manifold via $\chart^{-1}$.
For a given chart $\chart\from \Theta\to\R^{\dim(\Theta)}$ on a manifold
$\Theta$, and an abstract tensor $g$ at $\abstheta\in\Theta$, we denote
the coordinate expression of $g$ in chart $\chart$ by
$\Tang\chart(g)$ (specifying $\theta$ is not needed since an abstract
tensor $g$ includes its basepoint information).
Given a numerical tensor $g$ with coordinates expressed in
the chart $\chart$, we denote $\Tang_\abstheta\chart^{-1}(g)$ the
corresponding abstract tensor at $\abstheta\in \Theta$.

\begin{defi}[ (Online natural gradient in charts on a manifold)]
\label{def:absnatgrad}
Let $\Theta$ be a smooth manifold. For each $t\geq 1$, let
$\absp_t(y|\abstheta)$ be a probabilistic model on some variable $y$,
depending smoothly on $\abstheta\in \Theta$.

For each time $t\geq 1$, let $\chart_t\from \Theta\to\R^{\dim(\Theta)}$ be
a chart on $\Theta$.

The \emph{online natural gradient descent} for the observations
$y_t$, in the sequence of charts $\chart_t$, maintains an element
$\abstheta_t\in \Theta$ and a metric tensor $\absJ_t$ at $\abstheta_t$, defined
inductively by
% \begin{align}
% \numtheta &\gets \chart_t(\abstheta_{t-1}),\qquad \numJ\gets
% \inchart{\absJ_{t-1}}{\chart_t}
% \\
% \numJ_t &\gets (1-\gamma_t) \numJ + \gamma_t \,\E_{y\sim 
% p_{\abstheta_{t-1}}(y)} \left[\frac{\partial \ln p_{\chart_t^{-1}(\numtheta)}(y)}{\partial
% \numtheta}^{\otimes 2}\right]
% \\
% \numtheta_t &\gets \numtheta -\eta_t \, \numJ_t^{-1}\transp{\left(
% \frac{\partial \absloss_t(\chart_t^{-1}(\numtheta))}{\partial \numtheta}
% \right)}
% \\
% \abstheta_t&\gets \chart_t^{-1}(\numtheta_t),\qquad
% \absJ_t\gets (\Tang\chart_t)^{-1}_{\abstheta_t}(\numJ_t)
% \end{align}
% 
% Or equivalently
% \begin{align}
% \numJ_t &\gets \Tang\chart_t
% \left((1-\gamma_t) \absJ_{t-1} + \gamma_t \,\E_{y\sim 
% p_{\abstheta_{t-1}}(y)} \left[\frac{\partial \ln p_{\abstheta}(y)}{\partial
% \abstheta_{t-1}}^{\otimes 2}\right]\right)
% \\
% \numtheta_t &\gets \chart_t(\abstheta_{t-1}) -\eta_t \, \numJ_t^{-1}
% \d\Phi_t
% \left(
% \frac{\partial \absloss_t(\abstheta)}{\partial \abstheta_{t-1}}
% \right)
% \\
% \abstheta_t&\gets \chart_t^{-1}(\numtheta_t),\qquad
% \absJ_t\gets (\Tang\chart_t)^{-1}_{\abstheta_t}(\numJ_t)
% \end{align}
% 
% Or equivalently
\begin{align}
\label{eq:applychart}
\numtheta &\gets \chart_t(\abstheta_{t-1}),\qquad \numJ\gets
\Tang\chart_t(\absJ_{t-1})
\\
\label{eq:absnumJ}
\numJ_t &\gets (1-\gamma_t) \numJ +\gamma_t \Tang\chart_t
\left(\E_{y\sim 
\absp_t(y|\abstheta_{t-1})} \left[\frac{\partial \ln \absp_t(y|\abstheta)}{\partial
\abstheta_{t-1}}^{\otimes 2}\right]\right)
\\
\label{eq:absnumtheta}
\numtheta_t &\gets \numtheta +
\eta_t\,\numJ_t^{-1}
\Tang\chart_t
\transp{\left(
\frac{\partial \ln \absp_t(y_t|\abstheta)}{\partial \abstheta_{t-1}}
\right)}
\\
\label{eq:leavechart}
\abstheta_t&\gets \chart_t^{-1}(\numtheta_t),\qquad
\absJ_t\gets \Tang_{\abstheta_t}\chart_t^{-1}(\numJ_t)
\end{align}
with learning rate $\eta_t$ and Fisher matrix decay rate $\gamma_t$.
\end{defi}

If the chart $\Phi_t$ is constant in time, then this reduces to the
ordinary online natural gradient on $\theta_t$
(Lemma~\ref{lem:seminumnatgrad} below).
Indeed, if
$\chart_t=\chart_{t+1}$ then applying
$\chart_t^{-1}$ at the last step then applying $\chart_{t+1}$ in the next
step cancels out, so this amounts to just disregarding $\abstheta$ and
working on $\numtheta$.

\section{Kalman Filtering}
\label{sec:kalman}

One possible definition of the extended Kalman filter is as follows
\cite[\S15.1]{simon2006kalmanbook}. We
are trying to estimate the current state of a dynamical system $s_t$
whose evolution equation is known but whose precise value is unknown; at
each time step, we have access to a noisy measurement $y_t$ of a quantity
$\hat y_t=h(s_t)$ which depends on this state.

The Kalman filter
maintains an approximation of a Bayesian posterior on $s_t$ given the
observations $y_1,\ldots,y_t$. The posterior distribution after
$t$ observations is approximated by a Gaussian with mean $s_t$ and
covariance matrix $P_t$. (Indeed, Bayesian posteriors always tend to
Gaussians asymptotically under mild conditions, by the Bernstein--von
Mises theorem \cite{van2000asymptotic}.) The Kalman filter
prescribes a way to update $s_t$ and $P_t$ when new observations become
available.

The Kalman filter update is summarized in Definition~\ref{def:kalman}
below.  It is built to provide the \emph{exact} value of the Bayesian
posterior in the case of \emph{linear} dynamical systems with Gaussian
measurements and a Gaussian prior. In that sense, it is exact at first
order.

\begin{defi}[ (Extended Kalman filter)]
\label{def:kalman}
Consider a dynamical system
with state $s_{t}$, inputs $u_{t}$ and outputs $y_{t}$,
\begin{equation}
s_{t}=f(s_{t-1},u_{t})+\gaussian(0,Q_{t}),
\qquad
\hat y_{t}=h(s_{t},u_{t}),\qquad
y_{t}\sim \pobs(y|\hat y_{t})
\end{equation}
where $\pobs(\cdot|\hat y)$ denotes an exponential family with mean parameter $\hat y$
(e.g., $y=\gaussian(\hat y,R)$ with fixed covariance matrix $R$).

The \emph{extended Kalman filter} for this dynamical system
estimates the current state $s_t$ given observations $y_1,\ldots,y_t$ in
a Bayesian fashion.
At each time, the Bayesian posterior distribution of the state given
$y_1,\ldots,y_t$
is approximated by a Gaussian $\gaussian(s_t,P_t)$ so that $s_t$ is the
approximate maximum a posteriori, and $P_t$ is the approximate posterior
covariance matrix. (The prior is $\gaussian(s_0,P_0)$ at time $0$.)
Each time a new observation $y_{t}$ is available, these estimates are
updated as follows.

The transition step (before
observing $y_{t}$) is
\begin{align}
s_{t|{t-1}} & \gets f(s_{t-1},u_{t})
\\
\label{eq:KFdefF}
F_{t-1} & \gets \left.\frac{\partial f}{\partial s}\right|_{(s_{t-1},u_{t})}
\\
\label{eq:transP}
P_{t|{t-1}} & \gets F_{t-1} P_{t-1} \transp{F_{t-1}} + Q_{t}
\\
\hat y_{t} &\gets h(s_{t|{t-1}},u_{t})
\end{align}
and the observation step after observing $y_{t}$ is
\begin{align}
E_{t} &\gets \text{sufficient statistics}(y_{t})-\hat y_{t}
\\
R_{t} &\gets \Cov (\text{sufficient statistics}(y)|\hat y_{t})
\label{eq:KFdefR}
% \\
% E_{t} &\gets R_{t} \transp{\left(\frac{\partial \ln
% p(y_{t}|\hat y_{t})}{\partial \hat y_t}\right)}
\intertext{where the sufficient statistics are those of the exponential
family $\pobs$ (for a Gaussian model $y=\gaussian(\hat y,R)$ with known
$R$ these are just the error $E_{t}
= y_{t}-\hat y_{t}$ and the covariance matrix $R_t=R$)}
H_t &\gets \left.\frac{\partial h}{\partial
s}\right|_{(s_{t|{t-1}},u_{t})} \label{eq:defH}
\\
%S_{t} &\gets H_{t} P_{t|{t-1}} \transp{H_{t}}+R_{t}
%\\
%K_{t} &\gets P_{t|{t-1}}\transp{H_{t}}S_{t}^{-1}
K_{t} &\gets P_{t|{t-1}}\transp{H_{t}}\left(H_{t} P_{t|{t-1}}
\transp{H_{t}}+R_{t}\right)^{-1}\label{eq:KFdefK}
\\
P_{t} &\gets \left(\Id-K_{t}H_{t}\right)P_{t|{t-1}} \label{eq:KFPupdate}
\\
s_{t} &\gets s_{t|{t-1}} + K_{t}E_{t}\label{eq:KFsupdate}
\end{align}
\end{defi}

Defining the output noise via an exponential family $\pobs(y|\hat y)$
allows
for a straightforward
treatment of various output models, such as discrete outputs (by letting
$\hat y$ encode the probabilities of each class) or Gaussians
with unknown variance.  In the
Gaussian case with known variance our definition is fully standard.
However, for continuous variables with non-Gaussian output noise, the definition of $E_t$ and $R_t$ above
differs from the practice of modelling non-Gaussian noise via a nonlinear
function applied to Gaussian noise.\footnote{Non-Gaussian output noise is often modelled in
Kalman filtering via a continuous nonlinear function applied to a
Gaussian noise \cite[13.1]{simon2006kalmanbook}; this cannot easily represent discrete random variables.
Moreover, since the filter linearizes the function around the $0$ value
of the noise \cite[13.1]{simon2006kalmanbook}, in that approach the noise is still
implicitly Gaussian, though with a state-dependent
variance.\label{ft:nongaussian}}

\begin{defi}[ (Pure fading-memory Kalman filter)]
\label{def:Q}
The \emph{pure fading-memory Kalman filter} consists in taking the process
noise $Q_t$ proportional to $P_{t|t-1}$, so that the noise on the
dynamics of $s_t$ is modeled to be proportional to the current
uncertainty on $s_t$. Specifically, given a sequence of weights
$\alpha_t\geq 0$, we call \emph{pure fading-memory extended Kalman filter} the choice
\begin{equation}
Q_t=\alpha_t F_{t-1}P_{t-1}\transp{F_{t-1}}
\end{equation}
so that the transition equation \eqref{eq:transP} for $P$ becomes
\begin{equation}
\label{eq:fmemtransP}
P_{t|{t-1}} \gets (1+\alpha_t) F_{t-1} P_{t-1} \transp{F_{t-1}}
\end{equation}
and $Q_t$ is proportional to $P_{t|t-1}$.
\end{defi}

In the Bayesian interpretation of the Kalman filter, this amounts to
giving more weights to the likelihood of recent observations: the weight for
the likelihood of previous observations decreases by a factor
$1/(1+\alpha_t)$ at each step, hence the name ``fading memory''. This
prevents the filter from ultimately growing stale and corresponds to
larger learning rates for the natural gradient descent. The choice
$\alpha_t=0$, on the other hand, corresponds to a learning rate $1/t$ for
the natural gradient descent.

\section{Statement of the Correspondence}

\begin{enonce}{Notation for the dynamical system}
We consider a dynamical system
with state $s_{t}\in \R^{\dim(s)}$, inputs $u_{t}\in \R^{\dim(u)}$ and dynamics $f$, namely,
\begin{equation}
\label{eq:syst2}
s_{t}=f(s_{t-1},u_{t})
\end{equation}
where $f$ is a smooth function from $\R^{\dim(s)}\times \R^{\dim(u)}$ to
$\R^{\dim(s)}$. Predictions $\hat y_t\in \R^{\dim(\hat y)}$ are made on
observations $y_t\in \R^{\dim(y)}$ via
\begin{equation}
\hat y_t=h(s_t,u_t)
\end{equation}
where $h$ is a smooth function from $\R^{\dim(s)}\times \R^{\dim(u)}$ to
$\R^{\dim(\hat y)}$, and the observation model on $y_t$ is
\begin{equation}
y_t\sim \pobs(y_t|\hat y_t)
\end{equation}
where $\pobs$ is some exponential family with mean parameter $\hat y_t$ (such
as a Gaussian with mean $\hat y_t$ and known variance).
\end{enonce}

We refer to Appendix~\ref{sec:expfam} for a reminder on exponential
families.

\begin{enonce}{Notation for natural gradient on trajectories}
Given a dynamical system as above and a sequence of inputs $(u_t)_{t\geq 1}$, we denote by $\Traj$ the
set of trajectories of the dynamical system, i.e., the set of sequences
$\traj=(s_t)_{t\geq 0}$ such that $s_t=f(s_{t-1},u_t)$ for all $t\geq 1$.

We also define the chart $\Phi_t$ that parameterizes trajectories by their
state at time $t$:
\begin{equation}
\chart_t\from \traj \mapsto s_t
\end{equation}

Each trajectory $\traj=(s_t)_{t\geq 0}\in \Traj$ defines a probability
distribution on observations by setting
\begin{equation}
\label{eq:defptraj}
\absp_t(y_t|\traj)\deq \pobs(y_t|\hat y_t)=\pobs(y_t|h(\chart_t(\traj),u_t))
\end{equation}
where $\hat y_t=h(\chart_t(\traj),u_t)$ is the prediction made at time $t$ from trajectory $\traj$.
\end{enonce}

Thus, we can apply the general definition of online natural gradient
(Definition~\ref{def:absnatgrad}) to the models $\absp_t$ in the charts
$\chart_t$: this provides an online natural gradient descent on $\traj$
given the observations $y_t$.

Thus, given an initial trajectory $\traj^0$ (typically defined by its
initial state $s^0_0$), the online natural gradient descent produces an
estimated trajectory $\traj^t$ after observing $y_1,\ldots,y_t$. Let us
denote
\begin{equation}
s_t\deq \chart_t(\traj^t),\qquad s_{t|t-1}\deq \chart_t(\traj^{t-1})
\end{equation}
the estimated states at time $t$ of the trajectories $\traj^t$ and
$\traj^{t-1}$, respectively.
In the notation of Definition~\ref{def:absnatgrad}, the space $\Theta$ is
$\Traj$, the parameter $\abstheta_t$ is $\traj^t$, its expression
$\theta_t$ in the chart $\Phi_t$ is $s_t$, and the intermediate value
$\theta$ is $s_{t|t-1}$.

Our goal is to show that these satisfy the
same evolution equations as in the extended Kalman filter. Namely, we
will prove the following.

\begin{thm}[ (Extended Kalman filter as a natural gradient on
trajectories)]
\label{thm:kalnat}
Consider a dynamical system as above, an initial state $s_0$, a sequence of inputs
$(u_t)_{t\geq 1}$, and a sequence of observations $(y_t)_{t\geq 1}$.

Let $s_t$ be the state estimated at time $t$ by the pure fading-memory Kalman
filter (Defs.~\ref{def:kalman} and \ref{def:Q}) with observations $(y_t)$, initial state $s_0$ and initial covariance matrix $P_0$.

Let $\traj^t\in \Traj$ be the trajectory estimated after $t$ steps of the
online natural gradient for the model $\absp_t(y_t|\traj)$ in the chart
$\chart_t$ (Def.~\ref{def:absnatgrad}), initialized at the trajectory
$\traj^0$ starting at $s_0$ ($\traj^0_0=s_0$), and with initial Fisher matrix $J_0$.

Assume that the initializations and hyperparameters (learning rate,
fading memory rate) of the two algorithms are related via
\begin{align}
P_0=\eta_0 J_0^{-1},\qquad
\eta_t=\gamma_t,\qquad
\frac{1}{\eta_t}=\frac{1}{1+\alpha_t}\,\frac{1}{\eta_{t-1}}+1
\end{align}

Then for all $t\geq 0$, the trajectory $\traj^t$ passes through $s_t$ at time $t$:
\begin{equation}
\traj^t_t=s_t
\end{equation}
and the Kalman covariance and Fisher matrix satisfy
$P_t=\eta_t J_t^{-1}$.
\end{thm}

Let us give a few examples to illustrate the relation between
hyperparameters: with $\alpha_t=0$ (no process noise in the Kalman
filter), the natural gradient learning rates must satisfy
$1/\eta_t=1+1/\eta_{t-1}$, which is satisfied by
$\eta_t=1/(t+\mathrm{cst})$. This is the classical asymptotic rate for
parameter identification in statistical theory; on a dynamical system it
can be realized only in the absence of noise in the system.  On the other
hand, a constant $\alpha_t=\alpha>0$ corresponds to 
a constant gradient learning rate $\eta_t=\frac{\alpha}{1+\alpha}$ (and any
other choice of $\eta_0$ defines by induction a sequence $\eta_t$ that
tends to this value).

The extended Kalman filter appears as a natural gradient with the
particular natural gradient setting $\eta_t=\gamma_t$: namely, the Fisher
metric decay rate is equal to the natural gradient learning rate. This is
commented in \cite{natkal} for the case $f=\Id$; in short, the natural
gradient maintains a second-order approximation of the log-likelihood of
recent observations as a function of the parameter, and $\eta_t=\gamma_t$
corresponds to using the same decay rate for old observations in the
first-order and second-order terms.

Note that $\eta_t\to1$ when $\alpha\to \infty$: infinite noise on the
process corresponds to infinite forgetting of the past, and in that case
the Kalman filter jumps to the maximum likelihood estimator for the
latest observation (estimated at second order)
\cite[p.~212]{simon2006kalmanbook}. With the natural gradient,
a learning rate of $1$ corresponds to directly jumping to the minimum for
quadratic functions (learning rates above $1$ overshoot with the natural
gradient).

\section{Proof of Theorem~\ref{thm:kalnat}}

The proof considers the online natural gradient from
Definition~\ref{def:absnatgrad} on the trajectory space of a dynamical
system, and gradually makes all elements more explicit until we are left
with the extended Kalman filter.

For the proof, we shall assume that the function $f(\cdot,u_t)$ that maps
$s_{t-1}$ to $s_{t}$ is invertible; this guarantees that we can indeed
parameterize trajectories by their value $s_t$ at any time $t$. In
particular, the quantities $F_t$ in the extended Kalman filter are
invertible.  Without this assumption we would have to consider
equivalence classes of trajectories having the same value at time $t$, as
a function of time, which would make notation substantially heavier. This
is not really needed: in the end all terms $F_t^{-1}$ vanish from the
expressions, and the statement of Theorem~\ref{thm:kalnat} only involves
the value $s_t$ and Fisher matrix $J_t$ at time $t$, not the past
trajectory back-computed from $s_t$.

\subsection{Online Natural Gradient in Charts: Explicit Updates}

Here we consider the general setting of Definition~\ref{def:absnatgrad}: 
$\Theta$ is a smooth manifold; for each $t\geq 1$, 
$\absp_t(y|\abstheta)$ is a probabilistic model on some variable $y$,
depending smoothly on $\abstheta\in \Theta$;
for each time $t\geq 1$,  $\chart_t\from \Theta\to\R^{\dim(\Theta)}$ is
some chart on $\Theta$.

\begin{lem}
\label{lem:seminumnatgrad}
Denote
\begin{equation}
\nump_t(y|\numtheta)\deq \absp_t(y|\chart_t^{-1}(\numtheta))
\end{equation}
the expression of the probabilistic model in the chart $\chart_t$.

Then the updates \eqref{eq:absnumJ}--\eqref{eq:absnumtheta} for $\numJ_t$ and $\theta_t$ in the online natural
gradient are equivalent to 
\begin{align}
\label{eq:numJ}
\numJ_t &\gets (1-\gamma_t) \numJ +\gamma_t \,
\E_{y\sim 
\nump_t(y|\numtheta)} \left[\frac{\partial \ln \nump_t(y|\numtheta)}{\partial
\numtheta}^{\otimes 2}\right]
\\
\label{eq:numtheta}
\numtheta_t &\gets \numtheta +
\eta_t\,\numJ_t^{-1}
\transp{
\frac{\partial \ln \nump_t(y_t|\numtheta)}{\partial \numtheta}
}
\end{align}
where $\numtheta$ and $\numJ$ are as in \eqref{eq:applychart} above.

In particular, if the chart $\Phi_t$ is the same at all times, then the
abstract online natural gradient reduces to the usual online natural
gradient, because \eqref{eq:applychart} and \eqref{eq:leavechart} cancel
each other out.
\end{lem}

\begin{dem}
This follows by applying Lemma~\ref{lem:tensorinchart} from
Appendix~\ref{sec:geodiff} to the function $\ln
\absp_t(y|\abstheta)$.
\end{dem}

By studying the effect of a change of chart from applying
\eqref{eq:leavechart} at one step and then \eqref{eq:applychart} at the
next step, we are ready to obtain fully explicit expressions for the
online natural gradient that do not refer to manifold points $\abstheta$
or abstract tensors. The structure is closer to the extended Kalman
filter.

\begin{lem}
\label{lem:numnatgrad}
Denote
\begin{equation}
\psi_t\deq \chart_{t+1}\circ \chart_t^{-1},
\qquad
\Psi_t\deq \left.\frac{\partial \psi_t(\numtheta)
}{\partial\theta}\right|_{\numtheta=\numtheta_t}
\end{equation}
the change of chart from $t$ to $t+1$, and its derivative. Also denote
$
\nump_t(y|\numtheta)\deq \absp_t(y|\chart_t^{-1}(\numtheta))
$ as in Lemma~\ref{lem:seminumnatgrad} above.

Then the online natural gradient descent in the charts $\Phi_t$ is
equivalent to
\begin{align}
\label{eq:numthetatrans}
\numtheta &\gets \psi_{t-1}(\numtheta_{t-1})
\\
\label{eq:numJtrans}
\numJ&\gets \transp{(\Psi_{t-1}^{-1})}\,J_{t-1}\,\Psi_{t-1}^{-1}
\\
\label{eq:numJup}
\numJ_t &\gets (1-\gamma_t) \numJ +\gamma_t \,
\E_{y\sim 
\nump_t(y|\numtheta)} \left[\frac{\partial \ln \nump_t(y|\numtheta)}{\partial
\numtheta}^{\otimes 2}\right]
\\
\label{eq:numthetaup}
\numtheta_t &\gets \numtheta +
\eta_t\,\numJ_t^{-1}
\transp{
\frac{\partial \ln \nump_t(y_t|\numtheta)}{\partial \numtheta}
}
\end{align}
\end{lem}

\begin{dem}
Consider the effect of following Definition~\ref{def:absnatgrad}: we
apply \eqref{eq:leavechart} at one step and then \eqref{eq:applychart}
at the
next step. Namely, we go from chart $\chart_{t-1}$ to $\chart_t$.
The transformation rule \eqref{eq:numthetatrans} for $\numtheta$ is a
direct consequence of this.
Similarly, the transformation rule \eqref{eq:numJtrans} for $\numJ$ follows from the change
of coordinate formula for a $(0,2)$-tensor, given in
Lemma~\ref{lem:changeofchart} in Appendix~\ref{sec:geodiff}, when going
from chart $\chart_{t-1}$ to $\chart_t$. The rest is copied from
Lemma~\ref{lem:seminumnatgrad}.
\end{dem}

\subsection{Online Natural Gradient on Trajectories of a Dynamical
System}

We now specialize these results to the main situation considered in this
text, that of observations of a dynamical system.

Let us translate Lemma~\ref{lem:numnatgrad} in this setting. We first
need to explicit the function $
\psi_t= \chart_{t+1}\circ \chart_t^{-1}$ and its derivative $\Psi_t$.

\begin{lem}
\label{lem:psiisf}
In the setting above, for any time $t\geq 1$ and state $s\in \R^{\dim(s)}$ we have
\begin{equation}
\chart_t(\chart_{t-1}(s))=f(s,u_t)
\end{equation}
\end{lem}

\begin{dem}
Indeed, $\chart_{t-1}(s)$ maps $s$ to the trajectory $\traj\in \Traj$ whose state
at time $t-1$ is $s$. Then the state at time $t$ of $\traj$ is $f(s,u_t)$
by definition of $\Traj$.
\end{dem}

It is then immediate to translate Lemma~\ref{lem:numnatgrad} in this
setting. Remember that we apply this lemma to $\Theta=\Traj$, $\theta_t=s_t$ and
$\theta=s_{t|t-1}$ by definition.

\begin{cor}[ (Explicit form of the online natural gradient for a dynamical
system)]
\label{cor:explicitnatgrad}
The online natural gradient descent for the dynamical system above, in
the sequence of charts $\chart_t$, is
equivalent to
\begin{align}
s_{t|t-1} &\gets f(s_{t-1},u_t)
\\
F_{t-1}&\gets  \frac{\partial f(s_{t-1},u_t)
}{\partial s_{t-1}}
\\
\hat y_t&\gets h(s_{t|t-1},u_t)
\\
\label{eq:numJtrans2}
\numJ&\gets \transp{(F_{t-1}^{-1})}\,J_{t-1}\,F_{t-1}^{-1}
\\
\numJ_t &\gets (1-\gamma_t) \numJ +\gamma_t \,
\E_{y\sim 
\pobs(y|\hat y_t)} \left[
\frac{\partial \ln \pobs(y|\hat y_t)}{\partial s_{t|t-1}}^{\otimes 2}\right]
\\
s_t &\gets s_{t|t-1} +
\eta_t\,\numJ_t^{-1}
\transp{
\frac{\partial \ln \pobs(y_t|\hat y_t)}{\partial s_{t|t-1}}
}
\end{align}
where the last expressions depend on $s_{t|t-1}$ via $\hat
y_t=h(s_{t|t-1},u_t)$.
\end{cor}

\begin{dem}
By Lemma~\ref{lem:psiisf}, the function $\psi_{t-1}=\chart_t\circ
\chart_{t-1}^{-1}$ appearing in Lemma~\ref{lem:numnatgrad} is
$f(\cdot,u_t)$. Therefore, its derivative $\Psi_{t-1}$ at point
$\theta_{t-1}=s_{t-1}$ is
\begin{equation}
\Psi_{t-1}= \frac{\partial f(s_{t-1},u_t)
}{\partial s_{t-1}}
\end{equation}
This provides the updates for $s_{t|t-1}$ and for $J$ in the statement.

Next, the probability distribution $\nump_t(y|\theta)$ appearing in
Lemma~\ref{lem:numnatgrad} is $\absp_t(y|\chart_t^{-1}(\theta))$ by
definition. Here $\theta=s_{t|t-1}$. In
our situation, $\absp$ is defined by \eqref{eq:defptraj} namely
$\absp_t(y|\traj)=\pobs(y|h(\chart_t(\traj),u_t))$. Therefore, we obtain
\begin{align}
\nump_t(y|\theta)&=\absp_t(y|\chart_t^{-1}(s_{t|t-1}))
\\&=\pobs(y|h(\chart_t(\chart_t^{-1}(s_{t|t-1})),u_t))
\\&=\pobs(y|h(s_{t|t-1},u_t))
\\&=\pobs(y|\hat y_t)
\end{align}
and this ends the proof.
\end{dem}

\subsection{The Kalman State Update as a Gradient Step}

Here we recall some results from \cite{natkal} on the Kalman filter. These results
interpret the update step in the Kalman filter as a gradient descent step
preconditioned by the covariance matrix $P$, and make the relationship
with the Fisher information matrix of the observation model $\pobs$.

This relies on the output noise model $\pobs(y|\hat y)$ being an
exponential family. This is satisfied in the most common case,
when the model for $y$ is Gaussian with mean $\hat y$, but also for other
types of model, such as categorical outputs where $\hat y$ is the vector
of probabilities of all classes.

The following statement is Proposition~6 in \cite{natkal}.

\begin{prop}[ (Kalman filter as preconditioned gradient descent)]
\label{prop:Kalmanasgrad}
The update of the state $s$ in a Kalman filter can be seen as an online gradient
descent on data log-likelihood, with preconditioning matrix $P_t$. More
precisely, the update~\eqref{eq:KFsupdate} is equivalent to
\begin{equation}
s_t=s_{t|t-1}+P_t \transp{\left(\frac{\partial \ln
\pobs(y_t|\hat y_t)}{\partial s_{t|t-1}}\right)}
\end{equation}
where this expression depends on $s_{t|t-1}$ via $\hat
y_t=h(s_{t|t-1},u_t)$.
\end{prop}

The next proposition is known as the \emph{information filter} in the Kalman filter
literature, and states that the observation step for $P$ is additive when
considered on $P^{-1}$ (see \cite[(6.33)]{simon2006kalmanbook} or Lemma~9
in \cite{natkal})

\begin{lem}[ (Information filter)]
\label{prop:infofilter}
The update \eqref{eq:KFdefK}--\eqref{eq:KFPupdate} of $P_t$ in the
extended Kalman filter is equivalent to
\begin{equation}
P_t^{-1} \gets P_{t|t-1}^{-1} + \transp{H_t} R_t^{-1} H_t
\end{equation}
(assuming $P_{t|t-1}$ and $R_t$ are invertible).
% 
% In particular, for
% static dynamical systems ($f(s,u)=s$ and $Q_t=0$), the
% whole extended Kalman filter \eqref{eq:KFdefF}-\eqref{eq:KFsupdate}
% is equivalent to
% \begin{align}
% P_t^{-1} &\gets P_{t-1}^{-1}+\transp{H_t} R_t^{-1} H_t
% \\
% s_t &\gets s_{t-1} - P_t \transp{\left(\frac{\partial \ell_t(y_t)}{\partial
% s_{t-1}}\right)}
% \end{align}
\end{lem}

% \begin{dem}%[ of Lemma~\ref{prop:infofilter}]
% The first statement is well-known for Kalman filters
% \cite[(6.33)]{simon2006kalmanbook}. Indeed, expanding
% the definition of $K_t$ in the update \eqref{eq:KFPupdate} of $P_t$, we have 
% \begin{equation}
% P_t=P_{t|t-1}- P_{t|{t-1}}\transp{H_{t}}\left(H_{t} P_{t|{t-1}}
% \transp{H_{t}}+R_{t}\right)^{-1} H_t P_{t|t-1}
% \end{equation}
% but this is equal to 
% $(P_{t|t-1}^{-1}+\transp{H_t} R_t^{-1} H_t)^{-1}$
%  thanks to the Woodbury matrix identity.
% %$\left(A+UCV \right)^{-1} = A^{-1} -
% %A^{-1}U \left(C^{-1}+VA^{-1}U \right)^{-1} VA^{-1}$,
% 
% The second statement follows from Proposition~\ref{prop:Kalmanasgrad} and the
% fact that for $f(s,u)=s$, the transition step of the Kalman filter is
% just $s_{t|t-1}=s_{t-1}$ and $P_{t|t-1}=P_{t-1}$.
% \end{dem}

The next result (Lemma~10 from \cite{natkal}) states that after each observation, the Fisher information matrix of
the latest observation is added to $P^{-1}$.

\begin{lem}
\label{lem:hrh}
For
exponential families $\pobs(y|\hat y)$, the term $\transp{H_t} R_t^{-1} H_t$
appearing in Lemma~\ref{prop:infofilter} is
equal to the Fisher information matrix of $y$ with respect to the state
$s$,
\begin{equation}
\transp{H_t} R_t^{-1} H_t=\E_{y\sim \pobs(y|\hat y_t)}\left[
\frac{
\partial \ln \pobs(y|\hat y_t)}{\partial s_{t|t-1}}^{\otimes 2}
\right]
\end{equation}
where this expression depends on $s_{t|t-1}$ via $\hat
y_t=h(s_{t|t-1},u_t)$.
\end{lem}

By collecting these results into the definition of the Kalman filter, one
gets the following reformulation, which brings it closer to a natural
gradient.

\begin{cor}
\label{cor:kalmanasgrad}
The extended Kalman filter can be rewritten as
\begin{align}
s_{t|{t-1}} & \gets f(s_{t-1},u_{t})
\\
F_{t-1} & \gets \left.\frac{\partial f}{\partial s}\right|_{(s_{t-1},u_{t})}
\\
P_{t|{t-1}} & \gets F_{t-1} P_{t-1} \transp{F_{t-1}} + Q_{t}
\\
\hat y_{t} &\gets h(s_{t|{t-1}},u_{t})
\\
P_t^{-1} &\gets P_{t|t-1}^{-1} +
\E_{y\sim \pobs(y|\hat y_t)}\left[
\frac{
\partial \ln \pobs(y|\hat y_t)}{\partial s_{t|t-1}}^{\otimes 2}
\right]
\\
s_t&\gets s_{t|t-1}+P_t \transp{\left(\frac{\partial \ln
\pobs(y_t|\hat y_t)}{\partial s_{t|t-1}}\right)}
\end{align}
where the last expressions depend on $s_{t|t-1}$ via $\hat
y_t=h(s_{t|t-1},u_t)$.
\end{cor}

In the pure fading-memory case, the update for $P_{t|t-1}$ is
$P_{t|t-1}\gets (1+\alpha_t) F_{t-1} P_{t-1} \transp{F_{t-1}}$. In that
situation, comparing this rephrasing of the Kalman filter with the
explicit form of the natural gradient in
Corollary~\ref{cor:explicitnatgrad} makes it clear that $J_t$ is
proportional to the inverse of $P_t$. This is made precise in the
following statement.

\begin{prop}
The Kalman algorithm in Corollary~\ref{cor:kalmanasgrad} in the pure fading-memory
case, and the natural gradient algorithm in
Corollary~\ref{cor:explicitnatgrad}, are identical under the identification
\begin{equation}
P_t=\eta_t J_t^{-1}
\end{equation}
provided the hyperparameters $\eta_t$, $\gamma_t$ and $\alpha_t$ satisfy the following relations:
\begin{equation}
\eta_t=\gamma_t, \qquad
\frac{1}{\eta_t}=\frac{1}{1+\alpha_t}\,\frac{1}{\eta_{t-1}}+1
\end{equation}
\end{prop}

In particular, if these algorithms are initialized at the same point
(same $s_0$, and $P_0=\eta_0 J_0^{-1}$), they will remain identical at all
times.

\begin{dem}
Define $\tilde J_t\deq \eta_t P_t^{-1}$; we want to show that $\tilde
J_t$ follows the same evolution equation as $J_t$.

If this holds, then the update of $s_t$ will be identical in the two
algorithms, since one uses $P_t$ and the user uses $\eta_t J_t^{-1}$ to
precondition the gradient.

From Corollary~\ref{cor:kalmanasgrad} in the pure fading-memory
case we get
\begin{align}
\tilde J_t&= \eta_t P_t^{-1}
\\&= 
\eta_t P_{t|t-1}^{-1} +
\eta_t \E_{y\sim \pobs(y|\hat y_t)}\left[
\frac{
\partial \ln \pobs(y|\hat y_t)}{\partial s_{t|t-1}}^{\otimes 2}
\right]
\\&= 
\frac{\eta_t}{1+\alpha_t} \transp{(F_t^{-1})}P_{t-1}^{-1}F_t^{-1} +
\eta_t \E_{y\sim \pobs(y|\hat y_t)}\left[
\frac{
\partial \ln \pobs(y|\hat y_t)}{\partial s_{t|t-1}}^{\otimes 2}
\right]
\\&= 
\frac{\eta_t}{1+\alpha_t} \transp{(F_t^{-1})}\,\frac{\tilde
J_{t-1}}{\eta_{t-1}}F_t^{-1} +
\eta_t \E_{y\sim \pobs(y|\hat y_t)}\left[
\frac{
\partial \ln \pobs(y|\hat y_t)}{\partial s_{t|t-1}}^{\otimes 2}
\right]
\end{align}
while the full update for $J_t$ in Cor.~\ref{cor:explicitnatgrad} is
\begin{equation}
J_t=(1-\gamma_t) \transp{(F_t^{-1})}\,
J_{t-1}F_t^{-1} +
\gamma_t \E_{y\sim \pobs(y|\hat y_t)}\left[
\frac{
\partial \ln \pobs(y|\hat y_t)}{\partial s_{t|t-1}}^{\otimes 2}
\right]
\end{equation}

Thus, the two updates coincide if
\begin{equation}
\gamma_t=\eta_t, \qquad 1-\eta_t=\frac{\eta_t}{(1+\alpha_t)\eta_{t-1}}
\end{equation}
and in this case, if the algorithms are identical at time $t-1$ then they
will be identical at time $t$. This ends the proof of the proposition and
of Theorem~\ref{thm:kalnat}.
\end{dem}

% \begin{dem}
% Let us omit time indices for brevity. We have $\ds\frac{\partial
% \ell(y)}{\partial s}=\frac{\partial \ell(y)}{\partial \hat y}\frac{\partial
% \hat y}{\partial s}=\frac{\partial \ell(y)}{\partial \hat y}H$.
% Consequently, $\ds \E_y
% \left[\frac{\partial \ell(y)}{\partial s}^{\otimes 2}\right]=\transp{H}\,
% \E_y
% \left[\frac{\partial \ell(y)}{\partial \hat y}^{\otimes 2}\right] H$. 
% The middle term $\ds \E_y
% \left[\frac{\partial \ell(y)}{\partial \hat y}^{\otimes 2}\right]$ is the
% Fisher matrix of the random variable $y$ with respect to $\hat y$.
% 
% Now, for
% an exponential family $y\sim p(y|\hat y)$ in mean parameterization $\hat y$, the Fisher
% matrix with respect to $\hat y$ is equal to the inverse covariance matrix of
% the sufficient statistics of $y$ (Appendix, \eqref{eq:fishbarT}), that
% is, $R_t^{-1}$.
% \end{dem}

\section{Continuous-Time Case: the Kalman--Bucy Filter as a Natural
Gradient}
\label{sec:cont}

Consider now a continuous-time model of a dynamical system with state
$s$ and control or input $u_t$, with evolution equation
\begin{equation}
\dot s_t=f(s_t,u_t),\qquad
\end{equation}
and we want to learn the current state of the system from observations
$y_t$. As before, the observations are modeled via an observation
function $h(s_t,u_t)$ plus noise,
\begin{equation}
\label{eq:whitenoisemodel}
y_t=h(s_t,u_t)+W_t 
\end{equation}
where $W_t$ is a white noise process with known covariance matrix $R_t$.

The continuous-time analogue of the extended Kalman filter for this
situation is the extended Kalman--Bucy filter, which can be described
\cite{wikipediakalmanbucy} by the
two evolution equations (which mix the transition and the observation
steps of the discrete-time case)
\begin{align}
\dot s_t&=f(s_t,u_t)+K_t(y_t-h(s_t,u_t))
\\
\dot P_t&=F_tP_t+P_t\transp{F_t}-K_tH_tP_t+Q_t
\end{align}
where
\begin{equation}
F_t\deq \frac{\partial f(s_t,u_t)}{\partial s_t},\qquad H_t\deq
\frac{\partial h(s_t,u_t)}{\partial s_t},\qquad K_t\deq
P_t\transp{H_t}R_t^{-1}
\end{equation}
Here $Q_t$ is the covariance of the noise used to model the uncertainty
on the transitions of the system (for instance, if $f$ is not known
exactly), namely $\d s_t=f(s_t,u_t)\d t+\d B_t$
with $B_t$ a Brownian motion with covariance matrix $Q_t$.

As in the discrete case, we will work
with the \emph{pure fading-memory} variant of the extended Kalman--Bucy
filter, which assumes 
\begin{equation}
\label{eq:Qcont}
Q_t=\alpha_t P_t
\end{equation}
where $\alpha_t\geq 0$ is a hyperparameter. This choice of $Q_t$ is canonical
in the absence of further information on the system.

We will recover this filter fully in the course of proving
Theorem~\ref{thm:natkalbucy} below, by starting with an abstract
definition of the continuous-time online
natural gradient, and making it explicit until we end up with the
Kalman--Bucy filter.

The proof reveals a feature of the Kalman--Bucy filter: namely, to
properly define it in a manifold setting, a choice of covariant
derivative is needed to transfer the covariance matrix $P_t$ at the
current point $s_t$, to a new covariance matrix at $s_{t+\d t}$; in a
manifold this is a non-trivial operation (for the consequences for Kalman
filtering, see for instance the discussion
and Fig.~9 in the review \cite{barrau2018invariant}). This results from
the need to keep the algorithm online, and not recompute the Fisher
matrix of past observations when the parameter is updated.

On the other hand, the
evolution of the state $s_t$ does not depend explicitly on a covariant
derivative (or choice of chart), contrary to the discrete-time case:
in the discrete-time case, a change of chart influences the
update of $s_t$ only at second order in the learning rate, and this
disappears in continuous time because the learning rates become
infinitesimal.

\paragraph{Natural gradient in continuous time.}
The statistical learning viewpoint on this problem is as follows: Each
trajectory $\traj=(s_t)_{t\geq 0}$ defines a probability distribution on
generalized functions\footnote{We will call \emph{generalized functions} the
elements of the Wiener space, i.e., a functional space in which samples of the white
noise live. These can be seen, for instance, as random distributions
against which functions can be integrated.
Namely, the white noise definition states that if $w_t$ is sampled from a real-valued white noise with unit variance, then for every deterministic function $f$, the
integral $\int f(t) \,w_t \d t$ is Gaussian with variance $\int f(t)^2 \d
t$ \cite[Thm 4.1]{jazwinski_book}. 
Intuitively, the white noise $w_t$ takes values $(1/\sqrt{\d t})
\gaussian(0,1)$ in each infinitesimal interval of size $\d t$.
If $B_t=\int w_t\d t$ is the Brownian motion with derivative $w_t$, then $\int
f(t) \,w_t \d t$ is the same as the Itô integral $\int f(t)\d B_t$. For the
vector-valued case: if $w_t$ has covariance matrix $R_t$, then for each
vector-valued $f$ the integral $\int \transp{w_t} f(t) \d t$ is Gaussian
with variance $\int \transp{f(t)}R_t f(t)\d t$.} $\obs=(y_t)_{t\geq 0}$ in the Wiener space via
the model \eqref{eq:whitenoisemodel}.

The trajectories $\traj$ of the system may be parameterized via their initial state
$s_0$. Then the problem of estimating $\traj$ from the observations
$\obs$ becomes a standard statistical estimation problem of estimating
$s_0$, and methods such as the natural gradient may be applied to
optimize $s_0$ knowing the observations.

\begin{defi}[ (Observation model, instantaneous Fisher matrix,
instantaneous log-likelihood)]
\label{def:obsmodel}
We call \emph{observation model} parameterized by $\theta$ in some
manifold $\Theta$, 
a
probability distribution on generalized functions $\obs=(y_t)_{t\in
[0;T]}$ over $[0;T]$, which is absolutely continuous with density $p(\obs|\theta)$ with respect
to the Wiener measure.

The Fisher information matrix of this model over $[0;T]$ is
\begin{equation}
\label{eq:defcontJ}
J_{[0;T]}(\theta)\deq \E_{\obs \sim p(\obs|\theta)} \frac{\partial \ln
p(\obs|\theta)}{\partial \theta}^{\otimes 2}
\end{equation}
if this quantity exists.

We define the \emph{instantaneous} Fisher information matrix to
capture the amount of information brought by $y_t$ at instant $t$, as
\begin{equation}
\label{eq:jt}
j_t(\theta)\deq \frac{\d}{\d t} J_{[0;t]}(\theta)
\end{equation}
and the \emph{instantaneous log-likelihood} of $\obs=(y_t)_{t\in [0;T]}$,
which captures the likelihood of $y_t$ knowing the model, as
\begin{equation}
\ell_t(\theta)\deq \frac{\d}{\d t} \ln p(y_{[0:t]}|\theta)
\end{equation}
provided these derivatives exist.
\end{defi}

The instantaneous log-likelihood $\ell_t$ is the continuous-time analogue of the
log-likelihood of the current observation $y_t$ used to update the
parameter in
Definition~\ref{def:natgrad}. Intuitively $\ell_t$ is equal to $\ln
p(y_{[t;t+\d t]}|\theta,y_{[0;t)})$. We will formalize this intuition in
Proposition~\ref{prop:contll} and
Corollary~\ref{cor:conttheta} below: this corollary shows that for the model
$y_t=h(s_t,u_t)+W_t$, the gradient of this log-likelihood is given by the
error $y_t-h(s_t,u_t)$.

The instantaneous Fisher matrix $j_t$ is the continuous-time analogue of the Fisher information matrix on
a single observation $y\sim p(y|u_t,\theta)$ at time $t$ used in Definition~\ref{def:natgrad}.\footnote{Intuitively this is equal to
\begin{equation}
\d t \,\E_{\obs \sim p(\obs |\theta)}
\frac{\partial \ln p_t(y_t|\theta)}{\partial \theta}^{\otimes 2}
\end{equation}
which is formally closer to Def.~\ref{def:natgrad};
but this latter expression is not fully rigorous because 
samples $\obs\sim p(\obs|\theta)$ include white noise and thus have infinite values of $y_t$, which are
compensated by the $\d t$ factor. This is why we use the rigorous
expression \eqref{eq:contJ} instead.
% In the discrete setting, we could
% likewise have used $J_{[0:t]}-J_{[0:t-1]}$ which is an equivalent way of
% writing the Fisher information brought by the latest sample at time $t$
% in Def.~\ref{def:natgrad}.
}

The continuous-time analogue of the online natural gradient descent
(Def.~\ref{def:natgrad}) may be defined as follows.

\begin{defi}[ (Online natural gradient in continuous time)]
\label{def:contnatgrad}
Let $\obs=(y_t)_{t\geq 0}$ be a continuous function of time. Let $D$ be a
covariant derivative on the manifold $\Theta$.
Given an observation model as above,
we define the \emph{online natural gradient} for learning $\theta$ based
on the observations $\obs$, as the solution of
\begin{align}
\label{eq:contJ}
\frac{D J_t}{\d t} &=-\gamma_t J_t + \gamma_t \,j_t(\theta_t)
%\,\E_{y\sim p_t(y|\theta)} \frac{\partial \ln p_t(y|\theta)}{\partial \theta}^{\otimes 2}
\\
\label{eq:conttheta}
\dot \theta_t &= \eta_t J_t^{-1} \,\transp{\frac{\partial
\ell_t(\obs|\theta)}
{\partial \theta}}
\end{align}
initialized at some $\theta_0\in \Theta$ with some positive definite
metric tensor $J_0$.
\end{defi}

The term $-\gamma_t J_t$ in the equation introduces a decay
factor on $J$ as in the discrete case.

\paragraph{On the covariant derivative $D$ in the online natural
gradient.}
The covariant derivative $D$ is the continuous-time analogue of the choice of
charts at each time $t$ used in the discrete case. In the continuous-time
case, it is needed only for $J_t$, not for $\theta_t$. Indeed,
\eqref{eq:conttheta} is a well-defined ordinary differential
equation in the manifold $\Theta$, whose right-hand term is a tangent
vector at $\theta_t\in\Theta$ (see Lemma~\ref{lem:applymetric}). \footnote{We have assumed that the
observations $y_t$ are ordinary functions, not elements of the Wiener
space; in the latter case, \eqref{eq:conttheta} would become a stochastic
differential equation, requiring particular treatment to make it
parameterization-independent in the manifold.}

On the Kalman filter side of the correspondence, the need to introduce a
covariant derivative or coordinate system for $J$ corresponds to the fact that the
Kalman covariance matrix $P_{t|t-1}$ is translated from $s_{t|t-1}$ to
$s_t$ in the Kalman filter; this translation makes no sense in a
Riemannian manifold.

A canonical choice for $D$ would be the Levi-Civita
covariant derivative associated with the metric $J(\theta)$. However,
this does not result in a convenient algorithm. In the Kalman--Bucy
filter, $D$ turns out to be the covariant derivative associated with the
chart $s_t$ at time $t$; in particular, this $D$ is time-dependent.

The \emph{non-online} natural gradient would use
\begin{equation}
\dot \theta_t =\eta_t J(\theta_t)^{-1}\transp{\frac{\partial \ln
p_t(y_t|\theta)}{\partial \theta}}
\end{equation}
instead, which does not depend on a choice of covariant derivative. However, the
Fisher matrix $J(\theta)$ is an average over the time interval $[0;t]$
(from the statistical learning point of view, the Fisher matrix is an
expectation over inputs $u_t$):
using $J(\theta_t)$ would necessitate to recompute an
integral over the past for each new value of $\theta_t$.
Instead, the online version reuses values computed at previous
times instead of recomputing the full Fisher matrix $J(\theta_t)$ for new
values of $\theta_t$. This is why some way of transferring $J$ from
previous values of $\theta_t$ to the current one is needed.

Thus, the appearance of a covariant derivative in
Definition~\ref{def:contnatgrad} results from the need for a convenient
\emph{online} algorithm.

% TODO KEEP? Intuitively a lot of ``infinite constants'' are
% hidden by the use of the Wiener measure and the decomposition of
% $p(\obs|\theta)$ as a time integral of $p_t(y_t|\theta)$. For instance,
% due to the white noise, the values of a sample funciton $\obs$ under the distribution $p$ are infinite at every
% instant. In
% Appendix~\ref{sec:contintuition} we give an
% intuitive, non-rigorous account of all these infinite constants by taking
% the limit of a discrete-time system with very small time steps $\deltat$;
% in the end
% the infinite constants compensate.
% 
% However, let us insist that
% Definition~\ref{def:contnatgrad} makes sense as soon as the function
% $p_t$ and finite-valued observations $y_t$ are
% specified: the intuitive explanations of
% Appendix~\ref{sec:contintuition} are \emph{not} required to apply
% \eqref{eq:contJ}--\eqref{eq:conttheta} and to prove the results below.

\paragraph{The Kalman--Bucy filter as a natural gradient.} The
correspondence between the online natural gradient and the Kalman--Bucy
filter is expressed as follows.

\begin{thm}
\label{thm:natkalbucy}
Consider a continuous-time dynamical system
with state $s_{t}\in \R^{\dim(s)}$, inputs $u_{t}\in \R^{\dim(u)}$ and dynamics $f$, namely,
\begin{equation}
\label{eq:contsyst}
\dot s_{t}=f(s_{t-1},u_{t})
\end{equation}
where $f$ is a smooth function from $\R^{\dim(s)}\times \R^{\dim(u)}$ to
$\R^{\dim(s)}$. We assume that the solutions are regular on some time
interval $[0;T]$ for some open domain of initial conditions $s_0\in
\R^{\dim(s)}$.
Define the prediction model
\begin{equation}
\label{eq:obsmodel2}
y_t=h(s_t,u_t)+W_t
\end{equation}
where $h$ is a smooth function from $\R^{\dim(s)}\times \R^{\dim(u)}$ to
$\R^{\dim(y)}$, and $W_t$ is a white noise process with covariance matrix
$R_t$.

Let $\Theta$ be the set of trajectories of the system, parameterized by
their initial condition $\theta=s_0$. Thus, each $\theta\in \Theta$
defines a trajectory $s_t(\theta)$ and a probability distribution on observations $\obs=(y_t)_{t\in
[0;T]}$ via \eqref{eq:obsmodel2}. For each time $t\geq 0$, let $D^t$ be
the covariant derivative on $\Theta$ associated with the chart $\theta\mapsto
s_t(\theta)$ (namely, the covariant derivative $D^t$ of a tensor is equal
to its ordinary derivative when expressed in the chart $s_t$).

Let $(y_t)_{t\in
[0;T]}$ be a smooth series of observations. Let $\theta_t$ be the
trajectory at time $t$ inferred by the online natural gradient
(Def.~\ref{def:contnatgrad}) with observations $(y_t)$, where the
covariant derivative used at time $t$ is $D^t$.

Then the state $s_t(\theta_t)$ inferred by the natural gradient at time
$t$, is the same as the state $s_t$ inferred at time $t$ by the
Kalman--Bucy filter with pure fading memory ($Q_t=\alpha_t P_t$ in the
Kalman--Bucy equations), provided both
are initialized at the same state $s_0$, with Kalman--Bucy initial covariance
$P_0=\eta_0 J_0^{-1}$, and provided the hyperparameters are
related via
\begin{equation}
\gamma_t=\eta_t,\qquad 
\dot \eta_t=\alpha_t\eta_t-\eta_t^2
%\alpha_t=\frac{\dot \eta_t}{\eta_t}+\eta_t
\end{equation}

Moreover, the Kalman--Bucy posterior covariance is related to the
expression of the Fisher matrix $J_t$ in chart $s_t$ via
$P_t=\eta_t J_t^{-1}$.
\end{thm}

Let us comment once more on the hyperparameter settings. First, the extended
Kalman--Bucy filter (with fading memory $Q_t=\alpha_t P_t$) is recovered
as an online natural gradient with parameters $\gamma_t=\eta_t$ for the
same reasons as in the discrete case.

Second, the equation on $\eta_t$ is satisfied, for instance, if
$\eta_t=\alpha_t$ for all $t$. Other solutions exist: solutions are
better found by writing the equation on $1/\eta_t$ instead of $\eta_t$.
For instance, 
the full-memory, noiseless
case $\alpha_t=0$ corresponds to the learning rate
$\eta_t=1/(t+\mathrm{cst})$, as in the discrete case.

\section{Proofs for Continuous Time}

The proof proceeds by working out more and more explicit expressions for
the natural gradient, until we end up with the Kalman--Bucy filter.

We first compute an explicit form for the instantaneous log-likelihood
and its gradient, for the case of the model $y_t=h(s_t,u_t)+W_t$. This is
mostly a direct application of the Cameron--Martin theorem
\cite{cameronmartin1944}, and
relates the gradient of the instantaneous loss $\ell_t$ to the error
$y_t-h(s_t,u_t)$ at time $t$.

\begin{thm}[ (Cameron--Martin theorem)]
\label{thm:cameronmartin}
Let $h$ be a smooth real-valued function on $[0;T]$. Let $W_t$ be a
white noise on $[0;T]$. Let $\mathcal{W}$ be the Wiener measure (the law
of $W_t$ in the Wiener space). Let $\mathcal{W}_h$ be the Wiener measure
translated by $h$, namely, the distribution of $h+W$ in the Wiener space.
Then $\mathcal{W}_h$ is absolutely continuous with respect to
$\mathcal{W}$, and its density at a function $\obs=y(t)$ is
\begin{equation}
\label{eq:cameronmartin}
\frac{\d \mathcal{W}_h}{\d \mathcal{W}}(\obs)=\exp\left(\int_{[0;T]} y(t) h(t)\d
t-\frac12 \int_{[0;T]} h(t)^2\d t\right)
\end{equation}
\end{thm}

\begin{dem}
This is a rephrasing of Theorem~1 in \cite{cameronmartin1944}; the
statement given here can be found as Theorem~1.2 in
\cite{kuo75_gaussianmeasures} applied to the abstract Wiener space on $L^2([0;T])$
with variance $t=1$. 

At an informal level,
$y\deq h+W$ is a Gaussian centered at $h$ while the white noise
is a Gaussian centered at $0$, so informally the ratio of the probability
densities is the ratio of these two Gaussians,
\begin{equation}
\exp
\left(
-\frac12 \int_{t=0}^T (y(t)-h(t))^2\,\d t\right)\exp\left( 
\frac12 \int_{t=0}^T y(t)^2\,\d t\right)
\end{equation}
but rigorously, the quantity $\int_{t=0}^T y(t)^2\d t$
is infinite under the white noise distribution. However, this quantity
cancels out between the two parts of the expression, resulting in the
Cameron--Martin theorem and the expression~\eqref{eq:cameronmartin}.

More rigorously, let $E$ be a measurable set in the abstract
Wiener space over $L^2([0;T])$. Then $\mathcal{W}_h(E)=\mathcal{W}(E-h)$ by definition
of $\mathcal{W}_h$. Therefore, we can apply Theorem~1.2 in
\cite{kuo75_gaussianmeasures} to $-h$, which gives the result.

(Note that
we express everything over the white noise $W$ instead of the Brownian
motion 
$B=\int W$, so the norm we use for the Wiener space is indeed the $L^2$
norm instead of the square norm of derivatives as found for instance in
Theorem~1.1 of \cite{kuo75_gaussianmeasures}.)
\end{dem}

\begin{prop}[ (Instantaneous log-likelihood)]
\label{prop:contll}
Let $s_t(\theta)$ be a set of trajectories smoothly parameterized by
$\theta\in \Theta$. Consider the probability distribution on generalized
functions $\obs$ defined by \eqref{eq:whitenoisemodel}, namely,
$y_t=h(s_t(\theta,u_t))+W_t$ with $W_t$ a white noise with covariance
matrix $R_t$.

Then this probability distribution has a density $p(\obs|\theta)$ with
respect to the Wiener measure. For continuous functions $\obs$, this
density satisfies
\begin{equation}
\label{eq:obsmodel}
\ln p(\obs|\theta)=\int_{t=0}^T \ln p_t(y_t|\theta)\d t
\end{equation}
where
\begin{equation}
\label{eq:pt}
\ln p_t(y_t|\traj)\deq \transp{y_t}R_t^{-1}h(s_t,u_t)-\frac12
\transp{h(s_t,u_t)}R_t^{-1}h(s_t,u_t)
\end{equation}

Consequently, the instantaneous log-likelihood of this model is equal to
$\ell_t(\obs|\theta)=\ln p_t(y_t|\traj)$.
\end{prop}

(The probability density $p_t(y_t|\theta)$ does not sum to $1$ over $y_t$
because $p_t$ is not a probability but a probability density wrt the
Wiener measure; the Wiener measure contains the Gaussian factor on
$y_t$.)

\begin{dem}
The probability distribution on $\obs$ is, by definition, a white noise
centered at $h(s_t(\theta),u_t)$, with covariance matrix $R_t$. After
changing variables by $R_t^{-1/2}$, we can assume without loss of
generality that $R_t=\Id$. Then, by applying
Thm.~\ref{thm:cameronmartin} in the vector-valued case, we find that
the density of the law of $\obs$ with respect to the Wiener measure is
% $\mathcal{H}$ be the Hilbert space of functions $\obs \from[0;T] \to
% \R^{\dim(y)}$ equipped with the norm $\int_t \transp{y_t}R_t^{-1}y_t\d
% t$.
% The Cameron--Martin theorem TODO for the Wiener space over
% $\mathcal{H}$ states that the distribution of white noise centered at $h$
% has
% a well-defined density with respect to the Wiener measure (white noise
% centered at $0$),
% and that
% this density is
\begin{multline}
\label{eq:wienerdensity}
p(\obs|\traj)=\\
\exp
\left(
\int_{t=0}^T \transp{y_t} R_t^{-1}
h(s_t,u_t)\d t 
 -\frac12 \int_{t=0}^T \transp{h(s_t,u_t)} R_t^{-1}
h(s_t,u_t)\d t  \right)
\end{multline}
which proves the claim.
\end{dem}

This immediately provides an explicit form for the parameter update
\eqref{eq:conttheta} in the definition of the online natural gradient.

\begin{cor}[ (Gradient of instantaneous log-likelihood)]
\label{cor:conttheta}
Let $\traj=(s_t)_{t\in [0;T]}$ be a set of trajectories smoothly parameterized
by $\theta\in \Theta$.
Consider the associated observation model
\eqref{eq:whitenoisemodel} as above.

Then the natural gradient parameter update \eqref{eq:conttheta} for this
model satisfies
\begin{equation}
\label{eq:gradpt}
\frac{\partial \ell_t (\obs|\theta)}{\partial \theta}=
\transp{(y_t-h(s_t(\theta),u_t))}R_t^{-1}\,\frac{\partial h(s_t,u_t)}{\partial
\theta}
\end{equation}
\end{cor}

In particular, the gradient step will try to change the value of $h(s_t,u_t)$ to
reduce the error $y_t-h(s_t,u_t)$, as expected.

\begin{dem}
This is a direct consequence of \eqref{eq:pt}.
\end{dem}

We now turn to the expression for the Fisher matrix.

\begin{prop}[ (Instantaneous Fisher matrix)]
\label{prop:contfish}
Let $\traj=(s_t)_{t\in [0;T]}$ be a set of trajectories smoothly parameterized
by $\theta\in \Theta$, and
consider the observation model $y_t=h(s_t,u_t)+W_t$ with $W_t$ a
white noise with covariance $R_t$.
Denote
\begin{equation}
G_t\deq \frac{\partial s_t(\theta)}{\partial \theta},\qquad
H_t\deq \frac{\partial h(s_t,u_t)}{\partial s_t}
\end{equation}

Then the Fisher matrix \eqref{eq:defcontJ} for this model is
\begin{equation}
J_{[0:T]}(\theta)=
\int_0^T \transp{G_t}\transp{H_t}
R_t^{-1}H_tG_t \d t
\end{equation}
and in particular, the instantaneous Fisher matrix \eqref{eq:jt} is equal to
\begin{equation}
j_t(\theta)=\transp{G_t}\transp{H_t} R_t^{-1} H_tG_t
\end{equation}
\end{prop}

In particular, if the trajectories are parameterized by their state $s_t$
at time $t$ then $G_t=\Id$ (for this particular $t$), and this is
the continuous-time analogue of Lemma~\ref{prop:infofilter}.

\begin{lem}
\label{lem:whiteint}
Let $W_t$ be a vector-valued white noise on an interval $[0;T]$, with
covariance matrix $R_t$. Let $f(t)$ and $g(t)$ be two vector-valued
deterministic functions on $[0;T]$. Then
\begin{equation}
\E\left[
\left(
\int_0^T \transp{W_t}  f(t) \d t 
\right)
\left(
\int_0^T \transp{W_t}  g(t) \d t 
\right)
\right]
=
\int_0^T \transp{f(t)}R_t \,g(t)\d t
\end{equation}
(where the integrals are in the Wiener or Itô sense).
\end{lem}

\begin{dem}[ of Lemma~\ref{lem:whiteint}]
First, consider the case of a real-valued white noise $w_t$ with unit
variance. The integral $\int f(t)w_t\d t$ is equal to $\int f(t)\d B_t$
where $B_t$ is the Brownian motion whose derivative is $w_t$ (namely $\d
B_t=w_t \d t$). It is known
\cite[(4.23 for deterministic $f$ and $g$)]{jazwinski_book} that
\begin{equation}
\E \left[\left(\int_0^T f(t)\d B_t\right)\left(\int_0^t g(t)\d
B_t\right)\right]=\int_0^T f(t)g(t)\d t
\end{equation}
which gives the result for dimension $1$ and unit variance.

Now a vector-valued white noise with covariance matrix $R_t$ can
be written as $W_t=R_t^{1/2} \transp{(w^1_t,\ldots,w^n_t)}$ where the
$w^i_t$ are independent real-valued white noises with unit variance.
The result follows by
applying the above to $R_t^{1/2}f(t)$ and
$R_t^{1/2}g(t)$ and summing over components.
\end{dem}

\begin{dem}[ of Proposition~\ref{prop:contfish}]
The Fisher matrix for this model is, by \eqref{eq:defcontJ} and
\eqref{eq:obsmodel},
\begin{equation}
J_{[0;t]}=\E_{\obs\sim p(\obs|\theta)}\left[
\left(
\frac{\partial}{\partial \theta}\int_0^T \ln p_t(y_t|\theta)\d t
\right)^{\otimes 2}
\right]
\end{equation}
and the expression \eqref{eq:pt} for $\frac{\partial}{\partial \theta}
\ln p_t(y_t|\theta)$ yields
\begin{equation}
J_{[0;t]}=\E_{\obs\sim p(\obs|\theta)}\left[
\left(
\int_0^T
\transp{(y_t-h(s_t,u_t))}R_t^{-1}\,\frac{\partial h(s_t,u_t)}{\partial
\theta} \d t
\right)^{\otimes 2}
\right]
\end{equation}

Now, the model $p_t$ was derived from the observation model
\eqref{eq:whitenoisemodel}: under this model, $y_t=h(s_t,u_t)+W_t$ with
$W_t$ a white noise with covariance $R_t$. Therefore,
$y_t-h(s_t,u_t)=W_t$ and
\begin{equation}
J_{[0;t]}=\E_{\obs\sim p(\obs|\theta)}\left[
\left(
\int_0^T
\transp{W_t}\,R_t^{-1}\,\frac{\partial h(s_t,u_t)}{\partial
\theta} \d t
\right)^{\otimes 2}
\right]
\end{equation}

Now we can apply Lemma~\ref{lem:whiteint} to the components of the derivative with respect
to $\theta$, namely
\begin{equation}
f(t)= R_t^{-1}\frac{\partial h(s_t,u_t)}{\partial \theta_i}
\end{equation}
and
\begin{equation}
g(t)= R_t^{-1}\frac{\partial h(s_t,u_t)}{\partial \theta_j}
\end{equation}
and we find that the $(i,j)$ entry of the Fisher matrix is
\begin{equation}
\int_0^T \transp{\frac{\partial h(s_t,u_t)}{\partial
\theta_i}}R_t^{-1}\frac{\partial h(s_t,u_t)}{\partial \theta_j}\,\d t
\end{equation}
hence the result.
\end{dem}

By putting these two results together, we get a more explicit form of the
natural gradient for sets of trajectories. 

\begin{cor}
Let $\traj=(s_t)_{t\in [0;T]}$ be a set of trajectories smoothly parameterized
by $\theta\in \Theta$, and
consider the observation model $y_t=h(s_t,u_t)+W_t$ with $W_t$ a
white noise with covariance $R_t$ (namely, Def.~\ref{def:obsmodel} with
$p_t(y_t|\theta)$ given by \eqref{eq:pt}).
Denote
\begin{equation}
G_t\deq \frac{\partial s_t(\theta)}{\partial \theta},\qquad
H_t\deq \frac{\partial h(s_t,u_t)}{\partial s_t}
\end{equation}

Let $\obs=(y_t)_{t\in [0;T]}$ be a smooth function. Then the online
natural gradient (Def.~\ref{def:contnatgrad}) for this model with
observations $\obs$ satisfies
\begin{align}
\label{eq:numcontJ}
\frac{D^t  J_t}{\d t} &= -\gamma_t J_t + \gamma_t \,\transp{G_t}\transp{H_t} R_t^{-1} H_tG_t
\\
\label{eq:numconttheta}
\dot \theta_t&=\eta_t \,J_t^{-1}\transp{G_t}\transp{H_t}R_t^{-1}
(y_t-h(s_t(\theta_t),u_t))
\end{align}
where in these expressions, $G_t$ is evaluated at $\theta_t$ and $H_t$ at
$s_t(\theta_t)$, and where $D^t$ is the covariant derivative associated with the
chart $\theta\mapsto s_t(\theta)$.
\end{cor}

This result is still somewhat non-explicit due to the covariant
derivative $D^t$. This will disappear
by using $s_t$ rather than $\theta$
as the parameterization of the trajectories at each time; this is more consistent with
an algorithmic implementation at time $t$, and with the form of the
Kalman--Bucy filter.

Assume that $\theta\mapsto s_t(\theta)$ is indeed
a chart, namely, that $s_t$ is smooth and one-to-one on its domain with
smooth inverse. (This is the case under the assumptions of
Thm.~\ref{thm:natkalbucy}: then $s_t(\theta)$ is the solution of an
ordinary differential equation with initial condition $\theta=s_0$, and
if the function $f$ defining the equation is regular, then the mapping
from $s_{t_1}$ to $s_{t_2}$ is a diffeomorphism on its domain.)
This implies that $G_t=\partial s_t(\theta)/\partial
\theta$ is invertible.

Let $\Jt{t}{t_0}\deq \Tang s_{t_0} (J_t)$ be the expression of $J_t$ in chart $s_{t_0}$. Since
$J_t$ is a $(0,2)$-tensor at $\theta_t$ we have
\begin{equation}
\Jt{t}{t_0}\deq \Tang s_{t_0}(J_t
)=\transp{(G_{t_0}(\theta_t)^{-1})} J_t G_{t_0}(\theta_t)^{-1}
\end{equation}
by Lemma~\ref{lem:changeofchart}
(where we interpret this as a matrix expression by just viewing $\theta$ as
another chart).

We will be particularly interested in $\Jt{t}{t}$, which 
% \begin{equation}
% \Js_t\deq \transp{(G_t^{-1})}J_tG_t^{-1}
% \end{equation}
% assuming $G_t$ is invertible.
represents the Fisher matrix with respect to the current state $s_t$
rather than $\theta$. For this, we first have to study how $\Jt{t}{t_0}$
evolves when the reference chart $t_0$ changes. This works most finely
when the
trajectories parameterized by $\theta$ satisfy a differential equation
$\frac{\partial}{\partial t}s_t(\theta)=f_t(s_t(\theta))$ for some $f_t$,
as is the case in the Kalman--Bucy filter.

\begin{lem}[ (Time-varying charts)]
\label{lem:evolvingchart}
Let $\traj=(s_t)_{t\in [0;T]}$ be a set of trajectories smoothly parameterized
by $\theta\in \Theta$. Assume that there exists a function $f_t(s)$ such
that
\begin{equation}
\frac{\partial s_t(\theta)}{\partial t}=f_t(s_t(\theta))
\end{equation}
and set $F_t(s)\deq \frac{\partial f_t(s)}{\partial s}$.

Let $J$ be a $(0,2)$-tensor at some $\theta \in \Theta$. Then
the expression $\Jt{}{t_0}$ of $J$ in the chart $s_{t_0}$ evolves as
\begin{equation}
\frac{\d \Jt{}{t_0}}{\d
t_0}=-\transp{F_{t_0}}\,\Jt{}{t_0}-\Jt{}{t_0}F_{t_0}
\end{equation}
where $F_{t_0}$ is evaluated at $s_{t_0}(\theta)$.
\end{lem}

(Note that $\Jt{t}{t_0}$ is an expression in coordinates, so we can
take its ordinary derivative without needing covariant derivatives.)

This expression is the continuous-time analogue of
\eqref{eq:numJtrans2}: it shows that the transition update
$P_{t|t-1}\gets F_{t-1}P_{t-1}\transp{F_{t-1}}$ in the discrete-time extended
Kalman filter (and likewise the $F_tP_t+P_t\transp{F_t}$ term in the
Kalman--Bucy filter) is just a result of reexpressing the covariance
matrix in a coordinate system corresponding to the current state.

\begin{dem}
By the coordinate expression for a $(0,2)$-tensor in chart
$s_{t_0}(\theta)$, we have
\begin{equation}
\Jt{}{t_0}
=\transp{(G_{t_0}(\theta)^{-1})} J G_{t_0}(\theta)^{-1}
\end{equation}

So we are left with studying the derivative of $G_{t_0}(\theta)$.
Using
$\frac{\d}{\d t_0}G_t^{-1}=-G_{t_0}^{-1}(\frac{\d}{\d
t_0}G_{t_0})G_{t_0}^{-1}$ we find
\begin{equation}
\frac{\d \Jt{}{t_0}}{\d t_0}=-
\transp{(\dot G_{t_0}G_{t_0}^{-1})}\Jt{}{t_0}-\Jt{}{t_0}(\dot G_{t_0}G_{t_0}^{-1})
\end{equation}
where we have abbreviated $\dot G_{t_0}=\frac{\d}{\d
t_0}G_{t_0}$.

Now, the derivative of $G_t(\theta)$ with respect to $t$ satisfies
\begin{align}
\frac{\d}{\d t} G_t(\theta)
&=\frac{\partial}{\partial t}\frac{\partial}{\partial \theta} s_t(\theta)
=\frac{\partial}{\partial \theta}\frac{\partial}{\partial t} s_t(\theta)
\\&=\frac{\partial}{\partial \theta} f_t(s_t(\theta))
=\frac{\partial f_t(s_t(\theta))}{\partial s_t}\frac{\partial
s_t(\theta)}{\partial \theta}
\\&=F_t (s_t(\theta))G_t(\theta)
\end{align}
or more synthetically, $\dot G_t=F_t G_t$. This proves the claim.
\end{dem}

Then the online natural gradient rewrites as follows. Note the similarity
with the Kalman--Bucy filter for the state $s_t$, and for the Fisher
matrix $\Jt{t}{t}$ (which will ultimately be proportional to hte inverse
of $P_t$).

\begin{cor}[ (Explicit online natural gradient)]
Let $\traj=(s_t)_{t\in [0;T]}$ be a set of trajectories smoothly parameterized
by $\theta\in \Theta$. Assume that there exists a function $f_t(s)$ such
that
\begin{equation}
\frac{\partial s_t(\theta)}{\partial t}=f_t(s_t(\theta))
\end{equation}
Consider the observation model $y_t=h(s_t,u_t)+W_t$ with $W_t$ a
white noise with covariance $R_t$ (namely, Def.~\ref{def:obsmodel} with
$p_t(y_t|\theta)$ given by \eqref{eq:pt}).
Denote
\begin{equation}
F_t(s)\deq \frac{\partial f_t(s)}{\partial s},\qquad G_t\deq \frac{\partial s_t(\theta)}{\partial \theta},\qquad
H_t\deq \frac{\partial h(s_t,u_t)}{\partial s_t}
\end{equation}
and assume that $G_t$ is invertible.

Let $\obs=(y_t)_{t\in [0;T]}$ be a smooth function. Then the online
natural gradient (Def.~\ref{def:contnatgrad}) for this model with
observations $\obs$ is equivalent to
\begin{align}
\label{eq:numcontJ2}
\frac{\d}{\d t} \Jt{t}{t} &= -\transp{F_t}\Jt{t}{t}-\Jt{t}{t}F_t-\gamma_t \Jt{t}{t} + \gamma_t \,\transp{H_t} R_t^{-1} H_t
\\
\label{eq:numconttheta2}
\dot \theta_t&=\eta_t \,G_t^{-1}\Jt{t}{t}^{-1}\transp{H_t}R_t^{-1}
(y_t-h(s_t(\theta_t),u_t))
\end{align}
initialized with $\Jt{0}{0}\deq \transp{(G_0^{-1})}J_0G_0^{-1}$. In these
expressions, $F_t$ and $H_t$ are evaluated at $s_t(\theta_t)$ while $G_t$
is evaluated at $\theta_t$.

Moreover %$\Js_t$ is related to $J_t$ via $J_t=\transp{G_t}\Js_t G_t$, and 
the state $s_t(\theta_t)$ learned at
time $t$ satisfies
\begin{equation}
\label{eq:contstate}
\frac{\d}{\d t} s_t(\theta_t)=f_t(s_t(\theta_t))+
\eta_t \,\Jt{t}{t}^{-1}\transp{H_t}R_t^{-1}
(y_t-h(s_t(\theta_t),u_t))
\end{equation}
\end{cor}

\begin{dem}
Since $\Jt{t}{t}\deq \transp{(G_t^{-1})}J_tG_t^{-1}$, the relationship
\eqref{eq:numconttheta2} holds by direct substitution of the equation
evolution \eqref{eq:numconttheta} for $\theta_t$. 

By definition, $D^{t_0}$ is the covariant derivative associated with the
chart $s_{t_0}$ (Def.~\ref{def:chartconnection}). Its expression is obtained by going into the
chart $s_{t_0}$, taking ordinary derivatives, and going back, namely
(Def.~\ref{def:chartconnection}),
\begin{equation}
\frac{D^{t_0} J_t}{\d t}=\Tang s_{t_0}^{-1}\left(\frac{\d}{\d t}
\Jt{t}{t_0}\right)
\end{equation}

By definition of the online natural gradient with covariant derivative
$D=D^t$ at time $t$, one has, for $t_0=t$ 
\begin{equation}
\frac{D^{t_0} J_t}{\d t}=-\gamma_t J_t+\gamma_t j_t(\theta_t)
\end{equation}
Therefore, at time $t=t_0$, the expression above for $D^{t_0}$ yields
\begin{align}
\frac{\d}{\d t} \Jt{t}{t_0}&=\Tang s_{t_0}\left(
\frac{D^{t_0} J_t}{\d t}
\right)
\\&=\Tang s_{t_0}\left(
-\gamma_t J_t+\gamma_t j_t(\theta_t)
\right)
\\&=-\gamma_t \Jt{t}{t_0}+\gamma_t \transp{(G_{t_0}^{-1})}j_t G_{t_0}^{-1}
\end{align}
by definition of $\Jt{t}{t_0}$ and by the coordinate expression for
$(0,2)$-tensors in chart $s_{t_0}$. Here $G_{t_0}$ and $j_t$ are
evaluated at $\theta_t$.

Prop.~\ref{prop:contfish} provides the expression for the
instantaneous Fisher matrix $j_t$,
\begin{equation}
j_t(\theta)=\transp{G_t}\transp{H_t} R_t^{-1} H_tG_t
\end{equation}
so that for $t_0=t$, we have $\transp{(G_{t_0}^{-1})}j_t
G_{t_0}^{-1}=\transp{H_t} R_t^{-1} H_t$
and
\begin{align}
\frac{\d}{\d t} \Jt{t}{t_0}&=
-\gamma_t \Jt{t}{t_0}+\gamma_t \transp{H_t} R_t^{-1} H_t
\end{align}
at $t_0=t$.

Now we are interested in $\Jt{t}{t}$; to get its evolution equation we
must differentiate with respect to the two instances of $t$, one of which
captures the intrinsic change in $J$ and the other the change of chart:
\begin{align}
\frac{\d}{\d t} \Jt{t}{t}&=
\left.\frac{\d }{\d t}
\Jt{t}{t_0}\right|_{t_0=t}
+
\left.\frac{\d }{\d t_0}
\Jt{t}{t_0}\right|_{t_0=t}
\end{align}

We just computed $\left.\frac{\d }{\d t}
\Jt{t}{t_0}\right|_{t_0=t}$, and $\left.\frac{\d }{\d t_0}
\Jt{t}{t_0}\right|_{t_0=t}$ is provided by Lemma~\ref{lem:evolvingchart}.
This provides the full evolution equation \eqref{eq:numcontJ2} for
$\Jt{t}{t}$ in the statement.

To compute the evolution of the state $s_t(\theta_t)$ learned at
time $t$, let us decompose
\begin{align}
\frac{\d}{\d t} s_t(\theta_t)&=\left.\frac{\partial s_t(\theta)}{\partial
t}\right|_{\theta=\theta_t}+\frac{\partial s_t(\theta_t)}{\partial
\theta_t}\frac{\partial \theta_t}{\partial t}
\\&=f_t(s_t(\theta_t))+G_t\dot \theta_t
\end{align}
hence the result after substituting for $\dot \theta_t$.
\end{dem}

\begin{dem}[ of Theorem~\ref{thm:natkalbucy}]
First, by the assumptions of Theorem~\ref{thm:natkalbucy}, the trajectories $s_t$ satisfy the
evolution equation $\dot s_t=f(s_t,u_t)$, so we can apply the results
above with $f_t(s_t)\deq f(s_t,u_t)$, and the definition of $F_t$ is
consistent with the notation in the Kalman--Bucy filter.

Define
\begin{equation}
P_t\deq \eta_t \Jt{t}{t}^{-1}
\end{equation}
so that by definition, the evolution of the state \eqref{eq:contstate}
rewrites as
\begin{equation}
\frac{\d}{\d t} s_t(\theta_t)=f_t(s_t(\theta_t))+
P_t\transp{H_t}R_t^{-1}
(y_t-h(s_t(\theta_t),u_t))
\end{equation}
which is the state evolution equation in the Kalman--Bucy filter. Thus we are left
with checking the evolution equation for $P_t$.

We can compute the time derivative of $P_t$ via \eqref{eq:numcontJ2}. Using
$\frac{\partial}{\partial t} \Jt{t}{t}^{-1}=-\Jt{t}{t}^{-1}
(\frac{\partial}{\partial t}\Jt{t}{t})  \Jt{t}{t}^{-1}$, a direct computation
yields 
\begin{align}
\dot P_t&=\dot \eta_t \Jt{t}{t}^{-1}-\eta_t \Jt{t}{t}^{-1}
\left(
 -\transp{F_t}\Jt{t}{t}-\Jt{t}{t}F_t-\gamma_t \Jt{t}{t} + \gamma_t \,\transp{H_t} R_t^{-1} H_t
\right)
\Jt{t}{t}^{-1}
\\&=\frac{\dot \eta_t}{\eta_t} P_t+P_t \transp{F_t}+F_t P_t+\gamma_t
P_t-\frac{\gamma_t}{\eta_t}P_t \transp{H_t}R_t^{-1}H_t P_t
\end{align}

If $\gamma_t=\eta_t$, this coincides with the Kalman--Bucy evolution
equation for $P_t$ with process noise
\begin{equation}
Q_t=\left(\eta_t+\frac{\dot \eta_t}{\eta_t}\right)P_t
\end{equation}
which ends the proof.
\end{dem}

{\small

\appendix

\section{Appendix: Reminder on Exponential Families}
\label{sec:expfam}

\newcommand{\tsum}{{\textstyle \sum}}

An \emph{exponential family of probability distributions} on a variable
$x$ (discrete or continuous), with \emph{sufficient statistics}
$T_1(x),\ldots,T_K(x)$, is the
following family of distributions, parameterized by $\beta\in \R^K$:
\begin{equation}
p_\beta(x)=\frac{1}{Z(\beta)}\,\mathrm{e}^{\sum_k \beta_k T_k(x)}\,\lambda(\d x)
\end{equation}
where $Z(\beta)$ is a normalizing constant, and $\lambda(\d x)$ is any reference
measure on $x$, such as the Lebesgue measure or any discrete measure. The
family is obtained by varying the parameter $\beta\in \R^K$, called the
\emph{natural} or \emph{canonical} parameter. We will assume that the
$T_k$ are linearly independent as functions of $x$ (and linearly
independent from the constant function); this ensures that
different values of $\beta$ yield distinct distributions.

For instance, Bernoulli distributions are obtained with $\lambda$
the uniform measure on $x\in \{0,1\}$ and with a single sufficient
statistic $T(0)=0$, $T(1)=1$. Gaussian
distributions with a fixed variance are obtained with $\lambda(\d
x)$ the Gaussian distribution centered on $0$, and $T(x)=x$.

Another, often convenient parameterization of the same family is the
following: each value of $\beta$ gives rise to an average value $\bar T$
of the sufficient statistics, \begin{equation} \bar T_k\deq \E_{x\sim
p_\beta} T_k(x) \end{equation} For instance, for Gaussian distributions
with fixed variance, this is the mean, and for a Bernoulli variable this
is the probability to sample $1$.

Exponential families satisfy the identities
\begin{equation}
\label{eq:expder}
\frac{\partial \ln p_\beta(x)}{\partial \beta_k}=T_k(x)-\bar T_k,\qquad
\frac{\partial \ln Z}{\partial \beta_k}=\bar T_k
\end{equation}
by a simple computation \cite[(2.33)]{Amari2000book}.

These identities are useful to compute the Fisher matrix $J_\beta$ with
respect to the variable $\beta$, as follows \cite[(3.59)]{Amari2000book}:
\begin{align}
(J_\beta)_{ij} &\deq \E_{x\sim p_\beta} \left[
\frac{\partial \ln p_\beta(x)}{\partial \beta_i}
\frac{\partial \ln p_\beta(x)}{\partial \beta_j}
\right]
\\&= \E_{x\sim p_\beta} \left[ (T_i(x)-\bar T_i)(T_j(x)-\bar T_j)\right]
\\&= \Cov (T_i,T_j)
\label{eq:J=cov_comp}
\end{align}
or more synthetically
\begin{equation}
J_\beta=\Cov(T)
\end{equation}
where the covariance is under the law $p_\beta$. That is, for exponential
families the Fisher matrix
is the covariance matrix of the sufficient statistics. In particular it
can be estimated empirically, and is sometimes known algebraically.

In this work we need the Fisher matrix with respect to the mean parameter
$\bar T$,
\begin{equation}
(J_{\bar T})_{ij}=\E_{x\sim p_\beta}\left[
\frac{\partial \ln p_\beta(x)}{\partial \bar T_i}
\frac{\partial \ln p_\beta(x)}{\partial \bar T_j}
\right]
\end{equation}
By substituting $\frac{\partial \ln p(x)}{\partial \alpha}=\frac{\partial \ln
p(x)}{\partial \beta}\frac{\partial \beta}{\partial \alpha}$,
the Fisher matrices
$J_\alpha$ and $J_\beta$ with respect to parameterizations $\alpha$ and $\beta$ are
related to each other via
\begin{equation}
J_\alpha=\transp{\frac{\partial \beta}{\partial \alpha}}J_\beta\, \frac{\partial \beta}{\partial
\alpha}
\end{equation}
(consistently with the interpretation of the Fisher matrix as a
Riemannian metric and the behavior of metrics under change of coordinates
\cite[\S 2.3]{GHL87}).
So we need to compute $\partial{\bar T}/\partial \beta$. Using the
log-trick \begin{equation}
\partial \E_{x\sim p} f(x)=\E_{x\sim p} \left[f(x)\,\partial \ln
p(x)\right]
\end{equation}
together with \eqref{eq:expder}, we find
\begin{align}
\label{eq:expjac}
\frac{\partial \bar T_i}{\partial \beta_j}=
\frac{\partial \E T_i(x)}{\partial \beta_j}
=\E\left[ T_i(x)(T_j(x)-\bar T_j)
\right]=\E\left[ (T_i(x)-\bar T_i)(T_j(x)-\bar T_j)
\right]=(J_\beta)_{ij}
\end{align}
so that 
\begin{equation}
\frac{\partial \bar T}{\partial \beta}=J_\beta
\end{equation}
(see \cite[(3.32)]{Amari2000book},
where $\eta$ denotes the mean parameter) and consequently
\begin{equation}
\frac{\partial \beta}{\partial \bar T}=J_\beta^{-1}
\end{equation}
so that we find the Fisher matrix with respect to $\bar T$ to be
\begin{align}
J_{\bar T}&=\transp{\frac{\partial \beta}{\partial \bar T}}J_\beta\, \frac{\partial \beta}{\partial \bar T}
\\&=J_\beta^{-1} J_\beta J_\beta^{-1}
\\&=J_\beta^{-1}=\Cov(T)^{-1}
\label{eq:fishbarT}
\end{align}
that is, the Fisher matrix with respect to $\bar T$ is the inverse
covariance matrix of the sufficient statistics.

This gives rise to a simple formula for the natural gradient of
expectations with respect to the mean parameters. Denoting $\tilde
\nabla$ the natural gradient,
\begin{align}
\tilde\nabla_{\bar T}\, \E f(x) &\deq J_{\bar T}^{-1} \transp{\frac{\partial \E
f(x)}{\partial \bar T}}
\\&= J_{\bar T}^{-1} \,\transp{\frac{\partial \beta}{\partial \bar T}} \,\transp{\frac{\partial \E
f(x)}{\partial \beta}}
\\&=J_{\beta} J_\beta^{-1} \,\transp{\frac{\partial \E
f(x)}{\partial \beta}}
\\&=\transp{\frac{\partial \E
f(x)}{\partial \beta}}
\\&=\E\left[f(x) \,\frac{\partial \ln p_\beta(x)}{\partial \beta}
\right]
\\&=\E\left[f(x) (T(x)-\bar T)
\right]
\\&=\Cov(f,T)
\end{align}
which in particular, can be estimated empirically.

\section{Tensors and Charts for Manifolds}
\label{sec:geodiff}

We state without proof some classical results from differential geometry.

\begin{lem}
\label{lem:tensorinchart}
Let $\Theta$ be a smooth manifold, and let $\absloss\from \Theta\to \R$
be a real function on $\Theta$. Let $\abstheta\in \Theta$ and let $v$ be
the derivative of $\absloss$ at $\abstheta$, namely the cotangent vector
\begin{equation}
v=\frac{\partial \absloss(\abstheta)}{\partial \abstheta}
\end{equation}

Let $\chart\from \Theta\to \R^{\dim(\Theta)}$ be a chart on $\Theta$. Then
the expression of $v$ in the chart $\chart$ is
\begin{equation}
\Tang\chart(v)=\left.\frac{\partial \numloss(\numtheta)}{\partial
\numtheta}\right|_{\numtheta=\chart(\abstheta)}
\end{equation}
where
\begin{equation}
\numloss(\numtheta)\deq \absloss(\chart^{-1}(\numtheta))
\end{equation}
is the expression of $\absloss$ in the chart.

Similarly, the expression of $v^{\otimes 2}$ in the chart is
$\frac{\partial \numloss(\numtheta)}{\partial
\numtheta}^{\otimes 2}$.
\end{lem}

Remember that a $(0,2)$-tensor can be seen as a map sending a tangent
vector to a cotangent vector; therefore, if invertible, its inverse sends
cotangent vectors to tangent vectors.

\begin{lem}
\label{lem:applymetric}
Let $\absJ$ be an invertible $(0,2)$-tensor on a manifold $\Theta$, and
let $\mathfrak{v}$ be a cotangent vector at some $\abstheta\in \Theta$.
Let $\numJ$ and $v$ be respectively the matrix and row vector
representing $\absJ$ and $\mathfrak{v}$ in a chart. Then the expression
of $\absJ^{-1}\mathfrak{v}$ in the chart is $\numJ^{-1}\transp{v}$. 
\end{lem}

\begin{lem}
\label{lem:changeofchart}
Let $\Theta$ be a smooth manifold, and let $\absJ$ be a $(0,2)$-tensor at
some $\abstheta\in \Theta$. Let $\chart_1$, $\chart_2$ be two charts on
$\Theta$ and let $\psi\deq \chart_2\circ\chart_1^{-1}$ be the change of
chart.

Let $\numJ_1$ be the matrix representing $\absJ$ in chart $\chart_1$, and
likewise for $\numJ_2$. Then
\begin{equation}
\numJ_2=\transp{(\Psi^{-1})} \numJ_1\,\Psi^{-1}
\end{equation}
where
\begin{equation}
\Psi\deq \left.\frac{\partial \psi(\theta)}{\partial
\theta}\right|_{\theta=\chart_1(\abstheta)}
\end{equation}
\end{lem}

\begin{defi}[ (Covariant derivative associated with a chart)]
\label{def:chartconnection}
Let $\Theta$ be a smooth manifold and let $\chart\from \Theta\to
\R^{\dim(\Theta)}$ be a chart. The \emph{covariant derivative associated
with $\chart$} is the covariant derivative $D$ which coincides with the usual
derivative when expressed in chart $\chart$. Namely, for any curve
$(\theta_t)$ in $\Theta$ and any tensor $Z_t$ at $\theta_t$,
\begin{equation}
\frac{D Z_t}{\d t}\deq\Tang_{\theta_t} \chart^{-1}\left(
\frac{\d}{\d t} \Tang \chart(Z_t)
\right)
\end{equation}
\end{defi}

This is indeed a covariant derivative, whose Christoffel symbols in chart
$\chart$ are $0$.

}

\bibliographystyle{alpha}
\bibliography{kalnat}

\end{document}